\documentclass{amsart}
\usepackage{amssymb}

\title{A discrete model for Gell-Mann matrices}
\author{robert A. Wilson}
\date{First draft: 17th January 2024. This version 19th February 2024.}

\address{Queen Mary University of London}
\email{r.a.wilson@qmul.ac.uk}

\newcommand{\su}{\mathfrak{su}}
\newcommand{\FF}{\mathbb F}
\newcommand{\CC}{\mathbb C}
\newcommand{\RR}{\mathbb R}
\newcommand{\HH}{\mathbb H}
\newcommand{\rep}{\mathbf}
\begin{document}
\begin{abstract}
I propose a discrete model for the Gell-Mann matrices, which allows them to participate in discrete
symmetries of three generations of four types of elementary fermions, in addition to their usual role
in describing a continuous group $SU(3)$ of colour symmetries. This model sheds new light on the
mathematical (rather than physical) necessity for `mixing' between the various gauge groups $SU(3)$,
$SU(2)$ and $U(1)$ of the Standard Model. In particular it shows how the anti-Hermitian version of
Pauli matrices can act non-trivially on a unitary version of the Gell-Mann matrices, which 
leads to a non-trivial mixing between the weak and strong nuclear forces.

The unitary version of the Gell-Mann matrices
can in turn act non-trivially on
a quaternionic version of Dirac matrices, which leads to a non-trivial mixing between the strong force and 
the shape of spacetime defined by the Dirac matrices. Hence this model implies a mixing between the electro-weak-strong forces
on the one hand and gravity, as described by General Relativity, on the other. This mixing in turn implies the necessity for both
general relativistic corrections to the Standard Model of Particle Physics, and quantum corrections to General Relativity.
Contrary to general expectation, both types of corrections seem to be large enough to be tested experimentally.
\end{abstract}
\maketitle

\section{Introduction}
\label{intro}
\subsection{Context, aims and objectives}
The Gell-Mann matrices \cite{GellMann,Griffiths} are a particular choice of orthonormal basis for the (complex) Lie algebra $\su(3)$,
analogous to the Pauli matrices which form a basis for $\su(2)$. They are an essential part of the calculational
tools of quantum chromodynamics (QCD) as a description of the strong nuclear force \cite{QCD,WoitQFT}.
The analogy with Pauli matrices is not complete, however, since the Pauli matrices are unitary, but the
Gell-Mann matrices are not.
In particular, the Pauli matrices are
non-singular, so they can also be used to describe a \emph{finite} group of discrete (unitary) symmetries.
This is useful for describing the discrete symmetries of weak isospin, that distinguishes electrons from neutrinos.

 Specifically, the Pauli matrices generate a commuting product of the cyclic group $Z_4$ of order $4$ with the
 quaternion group $Q_8$ of order $8$, giving a finite analogue of the gauge group $U(1)_Y \times SU(2)_L$
of weak hypercharge and weak isospin. But the `mixing' of weak hypercharge and weak isospin to create
electric charge implies that $U(1)_{em}$ does not commute with $SU(2)_L$.  The only reasonable way to
accommodate this property in the finite groups is to take a cyclic group $Z_3$ in place of $U(1)_{em}$,
and let $Z_3$ act as automorphisms of $Q_8$. This leads to a group of order $24$, known as the binary tetrahedral group
(among many other names), which is a \emph{semi-direct} product $Q_8\rtimes Z_3$ of $Q_8$ and $Z_3$, rather than a direct product
$Q_8\times Z_3$.
Possible applications of this finite group to the modelling of elementary particles, and in particular of
electro-weak mixing, are explored in \cite{Yang,Frampton1,Frampton2,Frampton3,finite,tetrions}.

The Gell-Mann matrices, by contrast, are singular, so they do not generate a
finite group, and they cannot be used for discrete symmetries such as the three-generation symmetry of electrons.
This means that the three generations have to be added to the Standard Model `by hand', rather than arising automatically from the
formalism. The aim of this paper is to describe a method by which the Gell-Mann matrices can be given a discrete structure
that might avoid this problem. I also explore options for mixing the strong force with the electro-weak forces
by the process of using semi-direct products of the finite groups in place of the direct products of Lie groups.
It turns out that there is a unique possibility, arising from a particular 
three-dimensional complex reflection group \cite{ShepTodd,LehrerTaylor}
of order $648$.

The main objective of this paper 
is to explain the representation theory of this group of order $648$ in relation to the Standard Model.
An important part of this goal is to understand how the action of $Z_3$ on $Q_8$ relates to the Standard Model implementation of electro-weak mixing.
Similarly, the action of the Pauli matrices on the (modified) Gell-Mann matrices needs to be understood
in the context of weak-strong mixing. In all cases, the finite group can give only the combinatorial structure of this mixing, without
numerical values. The numerical values 
arise from choosing coordinates for the representations of the group.
It is at this point that the finite groups are embedded in unitary gauge groups, $U(1)$, $SU(2)$ and $SU(3)$ respectively,
and the structure of the Standard Model is recovered.

A secondary aim is to show how the (modified) Gell-Mann matrices can act on an extended (quaternionic) version of the Dirac matrices,
and thereby change the local shape of spacetime. We then investigate (at a qualitative level)
the possibility of using this action to develop a generally covariant model of quantum gravity, that appears to reduce to General Relativity \cite{GR3,GR1,GR2}
in the continuous limit,
if the electrons are ignored, and all matter is assumed to consist of neutrons.

\subsection{Pauli matrices}
The standard textbook \cite{Griffiths} Pauli matrices are 
\begin{align}
\sigma_x = \begin{pmatrix}0&1\cr 1&0\end{pmatrix}, \quad
\sigma_y = \begin{pmatrix}0&-i\cr i&0\end{pmatrix}, \quad
\sigma_z = \begin{pmatrix}1&0\cr 0&-1\end{pmatrix}
\end{align}
that is, the Hermitian matrices in the normal physics convention. 
These give rise to elements of $SU(2)$ via complex exponentiation as follows:
\begin{align}
&\exp(i\sigma_x\theta) = \begin{pmatrix}\cos\theta & i\sin\theta\cr i\sin\theta & \cos\theta\end{pmatrix}, \quad
\exp(i\sigma_y\theta) = \begin{pmatrix}\cos\theta & \sin\theta\cr -\sin\theta &\cos\theta\end{pmatrix},\cr
&\exp(i\sigma_z\theta) = \begin{pmatrix} e^{i\theta} & 0\cr 0&e^{-i\theta}\end{pmatrix}.
\end{align}
In addition there is a scalar matrix $i=\sigma_x\sigma_y\sigma_z$ which exponentiates to
\begin{align}
\exp(i\theta) = \begin{pmatrix}e^{i\theta} & 0\cr 0 & e^{i\theta}\end{pmatrix}
\end{align}
as elements of a scalar group $U(1)$. If this copy of $SU(2)$ is used for weak isospin, then the
corresponding copy of $U(1)$ is used for weak hypercharge.

Mathematicians prefer to put the multiplication by $i$ into the definition of the matrices,
instead of the exponentiation, so using
the anti-Hermitian matrices 
\begin{align}\label{KJI}
K:=\begin{pmatrix}0&i\cr i&0\end{pmatrix},\quad
J:=\begin{pmatrix}0&1\cr-1&0\end{pmatrix},\quad
I:=\begin{pmatrix}i&0\cr 0&-i\end{pmatrix}
\end{align}
in which the matrices $K,J,I$ behave like quaternions $IJ=-JI=K$ etc.

All the Pauli matrices, whether Hermitian or anti-Hermitian, are also unitary, as are the scalars. 
Therefore they generate finite subgroups of $U(2)$: in place of $U(1)_Y$ we have a cyclic group
$Z_4$ of order $4$ generated by the scalar $i$, and in place of $SU(2)_L$ we have a quaternion group $Q_8$ of order $8$,
generated by $i\sigma_x$, $i\sigma_y$ and $i\sigma_z$. 
In other words, there is a third interpretation of the Pauli matrices, as generators for (a unitary representation of) a finite group of order $16$, 
in addition to the
conventional interpretations as generators for unitary groups and Lie algebras. Such an interpretation is useful in cases
(such as weak isospin) in which the underlying symmetry is discrete.

\subsection{A $3$-dimensional analogue}
\label{3Dmats}
The Gell-Mann matrices in the normal physics convention are $3\times 3$ analogues of the Hermitian traceless Pauli matrices:
\begin{align}\label{GellMannmats}
&\lambda^1=\begin{pmatrix}0&1&0\cr 1&0&0\cr0&0&0\end{pmatrix},\quad
\lambda^2=\begin{pmatrix}0&-i&0\cr i&0&0\cr0&0&0\end{pmatrix},\quad
\lambda^3=\begin{pmatrix}1&0&0\cr0&-1&0\cr0&0&0\end{pmatrix}, \cr
&\lambda^4=\begin{pmatrix}0&0&1\cr 0&0&0\cr 1&0&0\end{pmatrix},\quad
\lambda^5=\begin{pmatrix}0&0&-i\cr 0&0&0\cr i&0&0\end{pmatrix},\quad
\lambda^6=\begin{pmatrix}0&0&0\cr0&0&1\cr0&1&0\end{pmatrix},\cr
&\lambda^7=\begin{pmatrix}0&0&0\cr 0&0&-i\cr0&i&0\end{pmatrix},\quad
\lambda^8=\frac1{\sqrt3}\begin{pmatrix}1&0&0\cr0&1&0\cr0&0&-2\end{pmatrix}.
\end{align}
Again, the usual mathematicians' convention is to multiply these matrices by $i$ to get anti-Hermitian matrices, which can then be exponentiated to
get unitary matrices generating $SU(3)$. In both cases, the Gell-Mann matrices form a basis for the $8$-space of complex $3\times 3$
 traceless matrices. 
 
 There are other possible bases, and it is suggested in \cite{octions} that a basis of \emph{real} matrices
might be more fundamental than Hermitian or anti-Hermitian matrices, for mathematical rather than physical reasons. 
But in none of these conventions are the matrices themselves unitary, so that the finite symmetries
one might hope to see (such as the generation symmetry for fundamental fermions) are not directly available in the standard formalism.
In fact it is easy to write down an orthonormal basis of unitary matrices, as I shall now demonstrate.

The natural analogue in three dimensions of the quaternion group (of order $2^3=8$) in two dimensions is the group $G_{27}$ of order $3^3=27$
generated by the matrices
\begin{align}\label{G27gens}
\begin{pmatrix}1&0&0\cr 0&\exp(2\pi i/3) & 0\cr 0&0&\exp(4\pi i/3)\end{pmatrix},\quad
\begin{pmatrix}0&1&0\cr 0&0&1\cr 1&0&0\end{pmatrix}.
\end{align}
Writing $v=\exp(2\pi i/3)$ and $w=v^2=\exp(4\pi i/3)$ for simplicity, the trace zero matrices form a complex $8$-space with the following
orthonormal basis:
\begin{align}\label{GellMannunitary}
&\begin{pmatrix}1&0&0\cr0&v&0\cr0&0&w\end{pmatrix},\quad
\begin{pmatrix}1&0&0\cr 0&w&0\cr 0&0&v\end{pmatrix},\quad
\begin{pmatrix}0&1&0\cr 0&0&1\cr 1&0&0\end{pmatrix},\quad
\begin{pmatrix}0&0&1\cr 1&0&0\cr 0&1&0\end{pmatrix},\cr
&\begin{pmatrix}0&1&0\cr 0&0&v\cr w&0&0\end{pmatrix},\quad
\begin{pmatrix}0&1&0\cr0&0&w\cr v&0&0\end{pmatrix},\quad
\begin{pmatrix}0&0&1\cr v&0&0\cr 0&w&0\end{pmatrix},\quad
\begin{pmatrix}0&0&1\cr w&0&0\cr 0&v&0\end{pmatrix}.
\end{align}

The first four matrices here are the generators (\ref{G27gens}) and their inverses, and the last four are the products of the first two with the next two.
The generators (of order $3$) do not commute with each other, and the commutators $x^{-1}y^{-1}xy$ are scalar matrices of order $3$.
In other words, there is a scalar group of order $3$ hidden in this group, consisting of scalars $1,v,w$, so that when
we exponentiate these matrices there is a threefold ambiguity in the overall phase, corresponding to the fact that the Lie group
$SU(3)$ has a centre consisting of scalars of order $3$, The above matrices are all unitary, and all have determinant $1$,
so that the group of order $27$ that they generate is a finite subgroup of $SU(3)$.

Moreover, these $8$ matrices span the same complex $8$-space as the $8$ Gell-Mann matrices, namely the space of all complex matrices with trace zero.
I therefore suggest that using these $8$ matrices in place of the $8$ Gell-Mann matrices 
makes essentially no difference to the Standard Model, as it is just a change of basis on the space of gluons.
But it has the huge advantage over the Gell-Mann matrices, of incorporating finite (triplet) symmetries in a natural way.
Indeed, it contains two 
\emph{independent} triplet symmetries, and may therefore be able to accommodate generation symmetries as
well as colour symmetries.

\section{Extensions}
\label{extensions}
\subsection{The binary tetrahedral group}
\label{theBTG}
The quaternion group has an automorphism of order $3$, that cycles the elements $I$, $J$ and $K$, and that can be represented
explicitly in the quaternion algebra by the matrix
\begin{align}
W:=(-1+I+J+K)/2 & = \frac12\begin{pmatrix}-1+i & 1+i \cr -1+i &-1-i\end{pmatrix}\cr
& = \frac{1+i}{2}\begin{pmatrix} i & 1\cr i & -1\end{pmatrix}.
\end{align}
This matrix extends the quaternion group $Q_8$ to a group of order $24$ that is known as the \emph{binary tetrahedral group}.
The above matrices give the `fundamental' or `natural' representation of this group inside $SU(2)$.

At the abstract group theory level, we have constructed a semi-direct product of $Q_8$ by $Z_3$, which we denote $Q_8\rtimes Z_3$. 
The subgroup $Q_8$ generated by the matrices $K$, $J$, $I$ in (\ref{KJI}) is normal, but the subgroup $Z_3$ generated by the above matrix $W$ is not.
In place of the relations 
\begin{align}
IZ=ZI, \quad JZ=ZJ
\end{align} that define the direct product, $Q_8\times Z_3 = \langle I,J\rangle \times \langle Z\rangle$, we have the relations 
\begin{align}
IW=WJ, \quad JW=WIJ
\end{align}
in the semidirect product $Q_8\rtimes Z_3 = \langle I,J\rangle \rtimes \langle W\rangle$.

The continuous analogue of $Q_8$ is $SU(2)$, and the
continuous analogue of $Z_3$ is $U(1)$, so that the effect of this semi-direct product is to create a
more complicated relationship between $SU(2)$ and $U(1)$ than the direct product $SU(2) \times U(1)$.
In physics language this relationship is a `mixing' of $SU(2)$ with $U(1)$ at the quantum level.
In the case of electro-weak mixing, $Z$ represents a generator for the group $U(1)_Y$ of weak hypercharge,
while $W$ represents a generator of the group $U(1)_{em}$ of (electromagnetic) charge, and $I, J, K$ represent the three components of
weak isospin, as generators for $SU(2)_L$. 

The usual `third component of weak isospin' is however more naturally associated with the imaginary part of $W$, that is  $(I+J+K)/2=W+1/2$,
corresponding to the usual equation $T_3=Q-Y_W/2$ relating charge $Q$ to weak hypercharge $Y_W$ and weak isospin $T_3$. Note that the
mathematics here indicates a preference for negating the weak hypercharge, compared to the usual convention.

\subsection{A ternary representation} 
In fact, the binary tetrahedral group can itself act as automorphisms of the group $G_{27}$ of order $27$ described in Section~\ref{3Dmats} above.
We may take matrices for the generators $I$ and $W$ as
\begin{align}\label{U2onGM}
I\mapsto-\frac{t}{3}\begin{pmatrix}1&1&1\cr 1&v& w\cr 1&w&v\end{pmatrix}, \quad W\mapsto \begin{pmatrix}v&0&0\cr 0&1&0\cr 0&0&1\end{pmatrix}
\end{align}where $t=v-w=\sqrt{-3}$, 
and check the defining relations
\begin{align}
I^2=-1, W^3=1, (IW)^3=1.
\end{align}
Generators for the quaternion subgroup $Q_8$ are then
\begin{align}\label{ternaryPauli}
I\mapsto-\frac{t}{3}\begin{pmatrix}1&1&1\cr 1&v& w\cr 1&w&v\end{pmatrix}, \quad
J\mapsto-\frac{t}{3}\begin{pmatrix}1&v&v\cr w&v& w\cr w&w&v\end{pmatrix},\quad
K\mapsto-\frac{t}{3}\begin{pmatrix}1&w&w\cr v&v& w\cr v&w&v\end{pmatrix}.
\end{align}

We must now check that these matrices act as automorphisms of $G_{27}$, which is easy in the case of $W$,
and for $I$ requires a few small 
calculations of this kind: 
\begin{align}
\begin{pmatrix}0&1&0\cr0&0&1\cr1&0&0\end{pmatrix}
\begin{pmatrix}1&1&1\cr 1&v& w\cr 1&w&v\end{pmatrix} & =
\begin{pmatrix} 1&v& w\cr 1&w&v\cr 1&1&1\end{pmatrix}\cr
& = \begin{pmatrix}1&1&1\cr 1&v& w\cr 1&w&v\end{pmatrix}
\begin{pmatrix}1&0&0\cr 0&v&0\cr 0&0&w\end{pmatrix}
\end{align}
which shows how the matrix for $I$ converts permutation matrices into diagonal matrices.
It follows that the matrices for $G_{27}$ and $Q_8$ together generate a group of order $216$ that combines finite analogues of $SU(3)$ and $SU(2)$ into a
single group. This is a semi-direct product $G_{27}\rtimes Q_8$, in which the subgroup $G_{27}$ is normal,
but the subgroup $Q_8$ is not. Since it is not a direct product, it does not give rise to  a direct product $SU(3) \times SU(2)$,
but instead `mixes' the two gauge groups.
If we compute the action of the matrices (\ref{ternaryPauli}) on the 3-vectors, we find that $(0,1,-1)$ is an eigenvector, with eigenvalue $1$.
Hence $Q_8$ embeds in a subgroup $SU(2)$ of $SU(3)$, which gives the required mixing.
 
The full group has order $2^3.3^4=648$, and is known to mathematicians as the triple cover of the Hessian group. 
It can be interpreted as a finite analogue of the complete gauge group
$U(1)\times SU(2)\times SU(3)$ of the Standard Model. 
But the finite group is not a direct product of three factors
in the way that the Lie group must be. It is an iterated semi-direct product
$G_{27}\rtimes Q_8\rtimes Z_3$. Thus the finite analogues of the three factors are mixed together in quite a complicated way,
which may have important consequences for the `mixing' of the different forces in the Standard Model.
We have already discussed electroweak mixing in the context of $Q_8\rtimes Z_3$, but here we have a much more
complicated weak-strong mixing in $G_{27}\rtimes Q_8$, as well as electro-strong mixing in $G_{27}\rtimes Z_3$.
Further details will be discussed below.

\subsection{A complex reflection group}
The group of order 648 is generated by the conjugates of $W$, i.e. the diagonal matrix in (\ref{U2onGM}), which has one eigenvalue $v$ and two eigenvalues $1$.
An invertible matrix with a single non-trivial eigenvalue is known as a \emph{reflection}. Real reflections have real eigenvalues, so the 
non-trivial eigenvalue is $-1$, but complex reflections can have any root of unity as an eigenvalue. Our group is therefore a
$3$-dimensional complex reflection group. There are $24$ reflections altogether, in $12$ `mirrors': 
\begin{align}\label{mirrors}
\begin{array}{cccc}
(t,0,0), & (1,1,1), & (v,1,1), & (w,1,1),\cr
(0,t,0), & (1,v,w), & (1,1,v), & (1,w,1),\cr
(0,0,t), & (1,v,w), & (1,v,1), & 1,1,w).
\end{array}
\end{align}
The two reflections in the mirror $r$ are obtained by putting  $\lambda=v$ or $\lambda=w$ in:
\begin{align}\label{reflform}
x \mapsto x - \frac{x.r}{r.r}(1-\lambda)r.
\end{align}

The four mirrors given in the top row of (\ref{mirrors}) generate one of the nine copies of the binary tetrahedral group, acting on the $2$-dimensional subspace 
perpendicular to $(0,1,-1)$.
In particular, this group is isomorphic to a $2$-dimensional complex reflection group. 
But we can also multiply the reflection by a scalar, to get two different copies of the binary tetrahedral group, which are not reflection groups,
and which will have different physical interpretations.
The images of the (anti-Hermitian) Pauli matrices can be written as products of two reflections of order $3$ (with opposite determinants).
This representation of $Q_8$ is equivalent to the sum of the trivial representation and the usual representation 
on a complex $2$-space or quaternion $1$-space. 
As a representation of the reflection group $Q_8\rtimes Z_3$ it splits as the trivial
representation plus a complex two-dimensional representation that is \emph{not} quaternionic.  One of the other two copies of $Q_8\rtimes Z_3$
acts as if on a complex number plus a quaternion, while the other does not. This appears to be a fundamental
mathematical reason for chirality of the weak interaction.

If we expand the complex $1$-spaces to real $2$-spaces, then the six complex scalar multiples of $12$ complex $3$-vectors become two real scalar multiples of
$36$ real $6$-vectors, and the corresponding $36$ real reflections generate the Weyl group of $E_6$. Hence the group 
$G_{27}\rtimes Q_8\rtimes Z_3$ is a subgroup of the Weyl group of $E_6$, and is in fact a maximal subgroup of the rotation part of the
Weyl group. It is possible that this fact may have led to the appearance of $E_6$ symmetry in certain Grand Unified Theories (GUTs) such as \cite{MD,Todorov}.
Unfortunately, however, the loss of the unitary structure of the $3$-dimensional representation appears to 
result in a loss of contact to real world physics. For this reason we are not considering models of $E_6$ type in this paper, but
sticking rigorously to models that can be justified by experiment.

\section{Representations and characters}
\label{representations}

\subsection{Continuous gauge groups}
So far we have only looked at the discrete, or combinatorial, structure of the combined group $G_{27}\rtimes Q_8\rtimes Z_3$, 
which gives a qualitative but not quantitative picture of the
underlying physics. In order to introduce measurable physical quantities such as mass, momentum and energy, and allocate them to elementary 
(or composite) particles we need to introduce real (or complex or quaternionic) representations, and use real (or complex or quaternionic) gauge groups
to choose coordinates for the representations.

There are three fundamental representations, described above, given by the embeddings of finite groups in Lie groups as follows:
\begin{align}
Z_3 & \subset U(1),\cr
Q_8\rtimes Z_3 & \subset Sp(1) \cong SU(2),\cr
G_{27}\rtimes Q_8 \rtimes Z_3 & \subset U(3).
\end{align}
These embeddings show how the gauge groups $U(1)$, $SU(2)$ and $SU(3)$ of the Standard Model
arise from the finite group. 
In each case, the real scalar factor has been taken out of the gauge group, since it only defines a unit of measurement, rather than a
genuine physical property. 

In the case of $U(1)$ acting on a complex $1$-dimensional representation, we are left therefore with one
physical parameter to measure, which can be represented as an angle, or phase, in the complex plane. In the case of $SU(2)$ acting on
a $1$-dimensional quaternionic representation, there are three such angles. Therefore by studying the representations of the binary tetrahedral
group $Q_8\rtimes Z_3$, we might expect to find four of the fundamental mixing angles of the Standard Model, including the electro-weak
mixing angle (Weinberg angle). A tentative identification of these four angles is obtained in \cite{finite}.

In the case of $U(3)$, there are five independent angles in 
the $3$-dimensional complex representation, bringing the total
number of mixing angles up to $9$,  
enough for the Weinberg angle plus four each in the CKM (Cabibbo--Kobayashi--Maskawa \cite{Cabibbo,KM})
and PMNS 
(Pontecorvo--Maki--Nakagawa--Sakata \cite{Pontecorvo,MNS}) matrices.
It is not obvious that all these $9$ mixing angles can in practice be found in these three fundamental representations, but 
I shall make the case below that at least $7$ of them arise naturally in this way.

\subsection{Clifford theory}
The standard method of constructing representations and characters of semi-direct products is Clifford theory \cite{Lux}, which is also used more generally
whenever a group has a (proper non-trivial) normal subgroup.
Since the group $G_{27}\rtimes Q_8 \rtimes Z_3$ has $4$ such normal subgroups, forming a single chain, 
we can work our way down the chain, from the largest normal subgroup
to the smallest, adding new representations at each stage. The method is standard, and the calculations straightforward,
so I do not give many details.

At the first stage, we have only the trivial quotient group, and the trivial representation.
At the second stage, we have a quotient $Z_3$, which adds a $2$-dimensional real representation, or equivalently a $1$-dimensional complex
representation and its complex conjugate. At the third stage, the quotient group is the alternating group $A_4$ on four letters, and the new representation is
the monomial representation of the group acting by conjugation on the three Pauli matrices, with signs attached. 

At the fourth stage, we have the fundamental pseudoreal representation as quaternions, with $Z_3$ as the corresponding `inertial quotient', so that the other representations
at this stage are obtained by tensoring with the representations of $Z_3$. This gives us a pair of complex conjugate two-dimensional representations,
or equivalently a real four-dimensional representation. Putting all these stages together gives us the representations of the binary tetrahedral group,
which are very well known and have been widely studied. The character table can be found in \cite[p.404]{JamesLiebeck} and is reproduced 
in Table~\ref{BTGchars}  for convenience,
in a shorthand form in which conjugacy classes that are scalar multiples of each other are grouped together.

\begin{table}\caption{\label{BTGchars}Character table of the binary tetrahedral group}
$$\begin{array}{l|cccc}
& \pm1 & K & \pm W & \pm W^2\cr\hline
\rep1a & 1 & 1 & 1 & 1\cr
\rep1b & 1 & 1 & v & w\cr 
\rep1c & 1  & 1 & w & v\cr 
\rep3 & 3 & -1 & 0 & 0\cr
\rep2a & \pm2 & 0 & \mp1 & \mp1\cr
\rep2b & \pm2 & 0 & \mp v & \mp w \cr 
\rep2c & \pm2 & 0 & \mp w & \mp v 
\cr\hline
\end{array}$$
\end{table}
Of particular significance here is that the three faithful representations are obtained by taking the fundamental representation $\rep2a$ and tensoring with
the representations $\rep1a$, $\rep1b$ and $\rep1c$ of $Z_3$. This type of pattern is a general feature of Clifford theory, and is a generalisation of the decomposition of representations
of direct products of groups as tensor products of representations of the factors. In particular,  we have a discrete version of the decomposition of
representations of $U(1)\times SU(2)$, but with the `mixing' of $U(1)$ with $SU(2)$ already incorporated into the mathematical structure.
A combinatorial version of the mixing occurs already at the group-theoretical level, but a quantitative estimate of the mixing angle can only be obtained
by looking in detail at the coordinates of the representations.

\subsection{The Hessian group}
Coming now to $G_{27}$, whose character table can be found in \cite[p.400]{JamesLiebeck}
we first take a $1$-dimensional complex representation, with inertial quotient $Z_3$, and induce up to the whole group
$G_{27}\rtimes Q_8\rtimes Z_3$, to get an $8$-dimensional
real representation, together with its tensor products with the representations of $Z_3$. 
The scalars in $G_{27}$ act trivially on these representations, so they are representations of a group
$(Z_3\times Z_3)\rtimes Q_8\rtimes Z_3$ of order $216$, which has been known as the Hessian group
since the 1870s.

The table of characters of the complex irreducible representations of 
the Hessian group is given in Table~\ref{PSU32chars}.
The vertical lines delineate normal subgroups, so that the first two columns comprise the normal subgroup $G_{27}/Z_3 \cong Z_3 \times Z_3$, 
generated by the modified Gell-Mann matrices
modulo scalars,
and the next two columns comprise the quaternion group $Q_8$ generated by the Pauli matrices. The horizontal lines separate cohorts of characters
for the various quotient groups. The top three rows contain the centralizer orders, the sizes of the conjugacy classes, and names for the classes, respectively.
The numerical part of the name gives the order of the group elements in that conjugacy class.

\begin{table}\caption{\label{PSU32chars}The character table of the Hessian group}
$$\begin{array}{r|rr|rr|rrrrrr|}
&216 & 27 & 24 & 4 & 18 & 9 & 6 & 18 & 9 & 6\cr
&1&8&9&54&12&24&36&12&24&36\cr
& 1A & 3A & 2A & 4A & 3B & 3C & 6A & 3B' & 3C' & 6A'\cr
\hline
\rep1a&1&1&1&1&1&1&1&1&1&1\cr
\rep1b&1&1&1&1&v&v&v&w&w&w\cr
\rep1c&1&1&1&1&w&w&w&v&v&v\cr
\hline
\rep3a&3&3&3& -1 & 0 &0&0&0&0&0\cr
\hline
\rep2a&2&2&-2&0&-1&-1&1&-1&-1&1\cr
\rep2b&2&2&-2&0 & -v &-v & v & -w & -w & w\cr
\rep2c&2&2&-2&0 & -w & -w & w & -v & -v & v\cr
\hline
\rep8a&8 & -1 & 0 & 0 & 2 & -1 & 0 & 2 & -1 & 0\cr
\rep8b&8&-1&0 & 0 & 2v&-v&0 & 2w&-w&0\cr
\rep8c&8&-1&0&0&2w&-w&0&2v&-v&0\cr\hline
\end{array}$$
\end{table}

For some purposes it is useful to replace a pair of complex conjugate characters by their real sum, and to fuse a conjugacy class with its inverse class.
This give us a simplified table (see Table~\ref{U32simple}).
\begin{table}\caption{\label{U32simple}Real characters of the Hessian group}
$$\begin{array}{c|rr|rr|rrr|}
&1&1&1&1&1&1&1\cr
U(1)=SO(2) &2&2&2&2&-1&-1&-1\cr\hline
SO(3)  & 3&3&3&-1&0&0&0\cr\hline
SU(2)  & 2&2&-2&0&-1&-1&1\cr
U(2) & 4 & 4 &-4 &0&1&1&-1\cr\hline
SU(3)& 8&-1&0&0&2&-1&0\cr
 U(3)& 16&-2&0&0&-2&1&0
\end{array}$$
\end{table}
The first column gives 
the name of the compact Lie group generated by the finite group, from which we can see a close
relationship to the gauge group of the Standard Model. The even-numbered rows contain the basic
building blocks $U(1)$, $SU(2)$ and $SU(3)$ respectively, and the other rows contribute to the
`mixing' between these three groups.

From this character table one can easily calculate the anti-symmetric square of $\rep8a$,
and we find
\begin{align}
\Lambda^2(\rep8a) = \rep2b+\rep2c+\rep8a+\rep8b+\rep8c.
\end{align}
In particular, there is a unique anti-symmetric (anti-commutative) algebra structure on $\rep8a$, namely the Lie algebra.
The symmetric square is similarly
\begin{align}
S^2(\rep8a) = \rep1a+\rep3a+\rep8a+\rep8a+\rep8b+\rep8c.
\end{align}
Hence there are two independent symmetric (commutative) products on $\rep8a$ as well.
One can therefore construct a wide variety of different algebra structures on $\rep8a$, some of which have been studied in the
physics literature \cite{Okubo,Okubonions}.

Indeed, the corresponding representations of the compact Lie group $SU(3)$ are as follows:
\begin{align}
\Lambda^2(\rep8) & = \rep8 +\rep{10} +\overline{\rep{10}}\cr
S^2(\rep8) & = \rep1 +\rep8 +\rep{27}
\end{align}
so that there is a (unique up to scalars) commutative algebra invariant under $SU(3)$, plus an infinite family of algebras that 
are neither commutative nor anti-commutative. 

\subsection{The faithful characters}

Finally we take a $3$-dimensional complex faithful representation of $G_{27}$, with inertial
quotient $Q_8\rtimes Z_3$, so that Clifford theory implies that we must now tensor with all the representations of the latter group. 
Let us arrange the conjugacy classes grouped together into scalar multiples. Then we only need to specify the character value on one of these classes,
and multiply by the appropriate scalars to get the others. In three cases an element is conjugate to its scalar multiples, so that the character values are all $0$,
while in the other seven they are not. 
The $14$ faithful complex characters are those given in Table~\ref{quarkchars} together with their complex
conjugates.

\begin{table}\caption{\label{quarkchars}Abbreviated table of characters for quarks}
$$\begin{array}{r|rr|rr|rrrrrr|}
\rep3b&3&0&-1&-1&t&0&t&-t&0&-t\cr
\rep3c&3&0&-1&-1&vt&0&vt&-wt&0&-wt\cr
\rep3d&3&0&-1&-1&wt&0&wt&-vt&0&-vt\cr\hline
\rep9&9&0&-3&1&0&0&0&0&0&0\cr\hline
\rep6a&6&0&2&0&-t&0&t&t&0&-t\cr
\rep6b&6&0&2&0&-vt&0&vt&wt&0&-wt\cr
\rep6c&6&0&2&0&-wt&0&wt&vt&0&-vt
\end{array}$$
\end{table}

It is not really necessary to write out this table in full, since, as noted above, the seven rows are obtained by multiplying the first row with the seven
irreducible characters of the binary tetrahedral group, that is the first seven characters of Table~\ref{PSU32chars}. If we take the latter to describe
properties of leptons, then Table~\ref{quarkchars} describes the corresponding properties of quarks.

Indeed, if we have the representations $\rep1b$, $\rep2a$ and $\rep3b$ then we can calculate all the other representations easily as tensor products,
subject only to the relations
\begin{align}
\rep2a \otimes \rep2a & = \rep1a +\rep3a\cr
\rep3b \otimes \overline{\rep3b} & = \rep1a + \rep8a.
\end{align}
Hence it is only necessary to gauge the representations $\rep1b$, $\rep2a$ and $\rep3b$, using the gauge groups $U(1)$, $SU(2)$ and $SU(3)$
respectively.

\section{Computing gauges}
\label{gauges}
\subsection{The electroweak gauge}
The $U(1)$ gauge group contains only one parameter, or `mixing angle', that identifies the directions of the real and imaginary parts of the representation.
This angle has been computed in \cite{finite} as $33.024^\circ$. The method of calculation was to superimpose the mass axis on an equilateral triangle
marked with the three generations of electrons on the vertices.
Hence it is expected that this mixing angle is one of the lepton mixing angles in the Standard Model.
However, it should be noted that this is an angle between the mass axis and one of the edges of the equilateral triangle, and it is
not clear \emph{a priori} which edge of the triangle to use, or in which direction. Of the six possible angles, only this one is
consistent with the experimental values of the mixing angles.
Hence we are not making a prediction here, but fitting the model to the experimental data, using a \emph{discrete} fitting parameter of $60^\circ$.

For convenience we summarise the calculation, which equates the physical mass ratio obtained from the three particles by swapping the electron and the muon
\begin{align}
\frac{\tau-e}{\tau-\mu} & = \frac{1776.86\pm.12-.51}{1776.86\pm.12 - 105.66}\cr
& = 1.062919\pm.000004
\end{align}
to the ratio of cosines of the angles between the mass axis and the two corresponding sides of the equilateral triangle:
\begin{align}
 \frac{\cos(60^\circ-\theta)}{\cos\theta} & = \cos60^\circ + \sin 60^\circ \tan\theta\cr
 & = \frac12 + \frac{\sqrt3}{2}\tan\theta.
\end{align}

The $SU(2)$ gauge group similarly contains exactly three independent parameters.
The angles calculated in \cite{finite} were $28.165516^\circ$, $66.7276^\circ$ and $2.337325^\circ$.
The method here was to identify discrete quaternions appropriate for four particles (the three generations of electron, plus a proton) that
are `fundamental' as far as electroweak interactions are concerned, and again superimpose a mass axis. 
There is again a discrete set of possible angles to measure, and again we are fitting the model to the experimental data
rather than making a specific prediction. But note that we are \emph{not} fitting any continuous parameters. The model we are considering here
has no continuous free parameters at all.

The calculations in this case are more complicated, but one of them has an obvious interpretation as a fundamental mass ratio:
\begin{align}
\cos\theta & = \frac{e}{n-p}\cr
 & \approx \frac{.510998950}{939.565420-.938.272088}\cr
 & \approx .395103
\end{align}
which gives the quoted value of $66.7276^\circ$. 

We should certainly expect one of
these angles to be the electroweak mixing angle $\theta_W$, also known as the Weinberg angle,
but it is not clear \emph{a priori} whether the other two are lepton-mixing (PMNS) 
or quark-mixing (CKM) angles. However, if we take into account that the embedding of $Q_8\rtimes Z_3$ into $SU(3)$ is described by one of the
representations $\rep1a+\rep2b$, $\rep1b+\rep2c$ or $\rep1c+\rep2a$, depending on which of the three copies of the binary tetrahedral group
we take, we should expect one `scalar' (i.e. CP-violating phase) and one generation-mixing angle.

It is therefore proposed in \cite{finite} to identify these angles, that are calculated from the directions of the mass axis in the two cases,
with four of the Standard Model mixing angles, as follows:
\begin{align}
\theta_{12}^{PMNS} & = {33.41^\circ}^{+.75^\circ}_{-.72^\circ},\cr
\theta_W & = 28.172^\circ \pm .021^\circ,\cr
\delta_{CP}^{CKM} & = 68.8^\circ\pm 4.5^\circ,\cr 
\theta_{23}^{CKM} & = 2.38^\circ\pm.06^\circ.
\end{align}

\subsection{The strong gauge}
The $SU(3)$ gauge group acts on a unitary $3$-space, which in principle requires $9$ `mixing angles' to define a unitary basis.
But we have already (conjecturally) identified $4$ of them from an embedding of $U(2)$ in $U(3)$. 
Hence there are just five more mixing angles to identify.
If the above conjectures are valid, then the remaining five
in the Standard Model are
\begin{align}
\theta_{12}^{CKM}&=13.04^\circ\pm.05^\circ,\cr
\theta_{13}^{CKM}&=.201^\circ\pm.011^\circ,\cr 
\theta_{13}^{PMNS}& = {8.54^\circ}^{+.11^\circ}_{-.12^\circ},\cr
 \theta_{23}^{PMNS}& = {49.1^\circ}^{+1.0^\circ}_{-1.3^\circ}, \cr
 \delta_{CP}^{PMNS}& = {197^\circ}^{+42^\circ}_{-25^\circ}.
\end{align}
In principle it should be possible to calculate these from an appropriate set of particle masses, but it is far from obvious how many masses we need,
or what the formulae are. Certainly we need at least six masses, which must presumably include all the quark masses, but it seems more likely that we need to
include the three electron masses and the proton mass as well, making ten in all.

A possible alternative is to use the baryon octet, which consists of eight particles with eight distinct masses, satisfying one known exact linear relation,
namely the Coleman--Glashow relation \cite{ColemanGlashow}:
\begin{align}\label{CGrelation}
p+\Sigma^-+\Xi^0 = n+\Sigma^++\Xi^-.
\end{align}
Current experimental values for the two sides of this mass equation, in units of MeV/$c^2$, are
\begin{align}
3450.58\pm.21 \sim 3450.64\pm.10
\end{align}
so that the agreement is well within current experimental uncertainty.
Indeed, it is possible to regard the Coleman--Glashow relation as picking out two triplets of particles that generalise the proton and the neutron,
analogous to the three generations of electron. Then it makes sense to arrange these nine particles in a $3\times 3$ array, as follows:
\begin{align}
\begin{array}{ccc}
e &  p & n\cr
\mu  & \Sigma^-& \Sigma^+\cr
\tau  & \Xi^0& \Xi^-
\end{array}
\end{align}
in such a way that the up/down/strange quark symmetry cycles the second and third columns in opposite directions. This symmetry belongs to the Hessian group,
thought of in its original form as the symmetry group of a $3\times 3$ affine plane.

\subsection{Identifying particles}
\label{particleids}
However it is done, we must start by identifying the particles we use with certain vectors or subspaces in $\CC^3$. Therefore we must begin by examining the
geometry of this $3$-dimensional space, using the complex reflection group as a guide.
The most important vectors are those that lie in the $12$ `mirrors'. The particular scalar multiple chosen is irrelevant for defining the reflection, so to begin with
at least we can ignore the scalar multiplication. The Hessian group then acts on this set of $12$ one-spaces in four triplets, as shown in (\ref{mirrors}).

The reflections in any mirror fix the mirrors in its column, and cycle the other three columns. The normal subgroup $G_{27}/Z_3$ of order $9$ contains four
subgroups of order $3$, each of which fixes the mirrors in one column, and cycles the mirrors in the other three columns (not always in the same direction!).
It is therefore tempting to identify these $12$ mirrors with the $12$ types of elementary fermions, and use the scalars for unobservable properties such as colours
of quarks and spin direction of electrons. However, the details of such an identification are tricky.

An alternative approach might be to fix a particular mirror, say $(t,0,0)$, so that the reflections rotate the mirror itself, to $(vt,0,0)$ and $(wt,0,0)$,
 fix the two perpendicular mirrors $(0,t,0)$ and $(0,0,t)$, and move the other $9$ mirrors. We can then perhaps identify the last $9$ mirrors as three generations of
 electrons and up and down quarks in some way, in order to get masses into the picture. Then we might use the mirror itself for a tenth mass such as a proton,
 or alternatively treat one or all of the three mirrors in the first column as neutrinos. The latter suggestion appears to mix the electron generation with the spin
 direction for a neutrino, which is exactly the kind of mathematical property that is needed both for the chirality of the weak interaction detected by the Wu experiment \cite{Wu}
 and for explaining neutrino oscillations \cite{SNO}.
 
 There is another set of $9$ one-spaces defined as the fixed spaces of the $9$ copies of $Q_8\rtimes Z_3$. These are the images of $(0,1,-1)$ under coordinate permutations and
 scalar multiplications on the coordinates. If we fix a particular copy of the group, then the $9$ one-spaces split as $1+8$, which suggests a splitting of the massless bosons
 in the Standard Model into one photon plus eight gluons. If we have identified the individual coordinates with neutrinos as suggested above, then there is a
 possible interpretation of gluons (and the photon) as bound states of a neutrino and an anti-neutrino. 
 Indeed, if we take the neutrinos as a basis for $\rep3b$ and the anti-neutrinos as a basis for $\overline{\rep3b}$, as above, then the bound states should lie in
 $\rep3b\otimes\overline{\rep3b} = \rep1a+\rep8a$.

Alternatively, if we use the Coleman--Glashow relation as a guide,
then we need three complex quantum numbers to represent the $8$ baryons. One possibility is 
\begin{align}
p=(v,-w,0),\quad & \Sigma^-=(0,v,-w),\quad \Xi^0=(-w,0,v),\cr
n=(-w,v,0),\quad & \Sigma^+=(0,-w,v), \quad \Xi^-=(v,0,-w).
\end{align}
Since these vectors span a real $5$-space, they are not quite enough to determine the mass vector uniquely. To get uniqueness, we must add in one of 
$\Lambda$ and $\Sigma^0$, if we want $(0,0,0)$ to have zero mass, or both, if we don't mind an offset from the origin.

\subsection{Calculating angles}
\label{moreangles}
It is not clear how best to do this, so let us begin with the six baryons above, and calculate the following:
\begin{align}
n-p=(1,-1,0) &= 1.293 \cr
\Sigma^--\Sigma^+=(0,-1,1) &= 8.08\pm.07\cr
\Xi^--\Xi^0 = (-1,0,1) & =6.5\pm.2 \cr
n+p=(t,t,0) & = 1877.837\cr
\Sigma^-+\Sigma^+ = (0,t,t) & = 2386.82\pm.08\cr
\Xi^-+\Xi^0 = (t,0,t) & = 2636.6\pm.2\cr
(t,t,t) & = 3450.6\pm.1\cr
(t,0,0) & = 1063.8\pm.1\cr
(0,t,0) & = 814.0\pm.1\cr
(0,0,t) & = 1572.8\pm.1
\end{align}

From these numbers we can calculate a finite number of angles, most of which 
have no obvious relationship to the Standard Model.
There is essentially only one angle to calculate from the $\Sigma$ baryons, that is given by
\begin{align}
(\Sigma^--\Sigma^+)/(\Sigma^-+\Sigma^+) & = .003385\pm.000029\cr
& = \sin(.194\pm.002)^\circ
\end{align}
which is consistent with the experimental value of $\theta_{13}$ in the CKM matrix, that is $(.201\pm.011)^\circ$, at about $0.6\sigma$.
This is a plausible interpretation, since it measures a mixing between the first-generation up/down quarks and the third-generation top/bottom quarks,
while keeping the strange and charm quarks fixed. Of course, this is in no sense a `prediction' of the model, but only a first attempt to fit the model
to experimental results. There are four more angles to identify.

If we assume that the Coleman--Glashow relation is exact, as I have done, then we have
\begin{align}
(n-p+ \Xi^--\Xi^0)/(\Sigma^+-\Sigma^-) & = -1\cr
& = \cos 180^\circ
\end{align}
giving a possible exact value for the $CP$-violating phase of the PMNS matrix.
If, on the other hand, we do not make this assumption, then the experimental value comes to around
\begin{align}
(7.8\pm.2)/(-8.08\pm.07) & = -.965\pm.025\cr
& \approx \cos({195^\circ}^{+5^\circ}_{-8^\circ})
\end{align}
which is in closer agreement with the experimental value, although this difference is not statistically significant, so we cannot draw any conclusions from it.

So far we have not used the fact that we are working in a complex $3$-space, but treated it merely as a real $6$-space.
The complex structure may well reduce the number of \emph{independent} mixing angles by up to $3$. The Coleman--Glashow relation removes one, so there may be
one or two more similar relations. If we exponentiate the mixing angles to get elements of $U(1)$, considered as eigenvalues of appropriate operators,
then multiplying eigenvalues together corresponds to adding angles, and adding three angles corresponds to calculating the determinant of a $3\times 3$ matrix.
For example, if we assume, as above, that the $CP$-violating phase of the PMNS matrix is exactly $180^\circ$, then the sum of the remaining three angles in
this matrix should represent a simple determinant:
\begin{align}
 {33.41^\circ}^{+.75^\circ}_{-.72^\circ}
+{8.54^\circ}^{+.11^\circ}_{-.12^\circ}
 +{49.1^\circ}^{+1.0^\circ}_{-1.3^\circ}
 = {91.05^\circ}^{+1.25^\circ}_{-1.5^\circ}.
\end{align}
It is therefore entirely possible that this angle is exactly $90^\circ$, representing a pure imaginary determinant, such as $t=\sqrt{-3}$.

At this stage we have three PMNS parameters, acted on by the subgroup $SO(3)$ of $SU(3)$ in its spin $1$ representation, and the other five Weinberg--CKM
parameters acted on in the spin $2$ representation. If there is another relation, it lies in this spin $2$ representation, and involves all five angles.
If we just add them up, we get around $112.6^\circ\pm4.5^\circ$, which could possibly be $120^\circ$, representing a complex number $v$ or $w$,
plausibly a product of two $3\times 3$ determinants. However the uncertainty in the $CP$-violating phase is so large that it is impossible to make a
sensible conjecture at this point. 
For example, we might equally well change the sign of the Cabibbo angle, giving $86.5^\circ\pm4.5^\circ$, and conjecture $90^\circ$.

Indeed, it seems more likely that no such relation exists between these angles. Altogether we have considered nine particle masses, that is three charged leptons
and six baryons, the latter forming a hexagram in the standard picture of the root system of $SU(3)$. 
The symmetries of the Standard Model make it clear that these nine particles consist of three `versions' of each of the particles of `normal' matter, that is the electron, proton and neutron.
The electron comes in three generations, that is the electron, muon and tau particle, while the proton comes in a triple with the $\Sigma^-$ and $\Xi^0$ baryons,
and the neutron is grouped together with the $\Sigma^+$ and $\Xi^-$ baryons.

These nine particles appear to form a basis for the adjoint
representation of $SU(3)$, which splits as $1+3+5$ for the real subgroup $SO(3)$. Here the scalar represents the determinant, and apparently gives us a single relation,
saying the determinant is pure imaginary, as we would expect from the pure imaginary mass term in the Dirac equation, if we have an odd number of masses.
The other two components should give us two more relations, which would appear to be the Coleman--Glashow relation and the equation 
\begin{align}
e+\mu+\tau+3p = 5n
\end{align}
that was pointed out in \cite{remarks}.
Therefore we should not expect any more universal mass relations among these nine particles, which makes it unlikely that there are any more
relations between the mixing angles either.

\section{Geometry}
\label{geometry}
\subsection{Finite fields} 
The underlying mathematical structure of the model is best elucidated by regarding $1,v,w$ not as complex numbers, but as abstract elements of the field of order $4$.
This is a field of characteristic $2$, so that $2=0$, $-1=1$ and $3=1$. The equation $1+v+w=0$ still holds, as in the complex numbers, so that $1+v=w$ and
$1+w=v$, as well as $v+w=1$. The multiplicative rules $v^2=w$, $w^2=v$ and $vw=1$ also hold just as before, and the involution that swaps $v$ with $w$
is the appropriate analogue of complex conjugation. Denote this field by $\FF_4$,
and denote the $n$-dimensional unitary group by $U(n,\FF_4)$.
Then it is easy to see that $U(1,\FF_4)$ is the cyclic group of order $3$, consisting of all non-zero $1\times 1$ matrices. 
The whole group of order $648$ constructed above is then $U(3,\FF_4)$. 

But the quaternion group and the binary tetrahedral group do not have obvious names in terms of $\FF_4$.
The reason for this is that the Pauli matrices are more naturally written over the field $\FF_9$ of order $9$, whose non-zero elements are $\pm1$, $\pm i$ and $\pm1\pm i$.
We can then see that the binary tetrahedral group is isomorphic to $SU(2,\FF_9)$. From this point of view, therefore, the reason for the complicated mixing between
$SU(2)$ and $SU(3)$ in the Standard Model is ultimately a complicated relationship between the finite fields of orders $4$ and $9$.
The particular feature that makes the finite groups behave completely differently from Lie groups is the chain of quotient maps 
\begin{align}
U(3,\FF_4) \rightarrow SU(2,\FF_9) \rightarrow U(1,\FF_4)
\end{align} 
This provides a unification of $U(3)$ and $SU(2)$ that is specific to these finite fields, and occurs nowhere else in mathematics.
In particular, such unification cannot be achieved with Lie groups over the complex numbers.

The $3$-dimensional space over $\FF_4$ on which $U(3,\FF_4)$ acts contains $4^3-1=63$ vectors, each in a set of $3$ scalar multiples, and therefore forming
$63/3=21$ one-dimensional subspaces. These $21$ `points' of the geometry divide into $12$ `odd' (non-singular) and $9$ `even' (singular) points. Each point is
perpendicular to a `line'. The $12$ non-singular lines each contain $2$ non-singular and $3$ singular points, while the $9$ singular lines each contain $1$
singular point and $4$ non-singular points.

The points can be arranged as
\begin{align}\label{points}
\begin{array}{|cccc|}
100 & 111 & v11 & w11\cr
010 & 1vw & 11v & 1w1\cr
001 & 1wv & 1v1 & 11w
\end{array}
\qquad
\begin{array}{|ccc|}
011 & 101 & 110\cr
0vw & w0v & vw0\cr
0wv & v0w & wv0
\end{array}
\end{align}
to exhibit the symmetry. 
These are just the reductions to $\FF_4$ of the $12$ `mirrors' (corresponding to $Z_3$ subgroups) and $9$ copies of the Pauli matrices
(corresponding to $Q_8$ subgroups), so that all of the above discussion of particles (Section~\ref{particleids}) applies equally to the finite geometry.

\subsection{Interpretations and speculations}
In particular,
we can interpret the $12$ non-singular points as the fundamental fermions, arranged into four columns of three generations (flavours) each,
and it would seem to make sense to interpret the first two columns as leptons, and the last two as quarks, although the mathematics allows us an arbitrary choice here.
It is noticeable, however, that in this interpretation
there is no single `generation symmetry' that acts on all four columns simultaneously. Instead, there is a group of order $9$ containing
elements of order $3$ that act on any three of the four columns.

This is physically reasonable, since it is known that quark generations behave differently from lepton generations in the weak interaction, so that at least
two different generation symmetries are required in order to describe what is experimentally observed. We can choose the two generation symmetries in
a variety of ways, for example an electron generation symmetry that acts on electrons but not on neutrinos, and a neutrino generation symmetry
that acts on neutrinos but not on electrons. Such a choice shows explicitly that the neutrinos can change flavour \emph{independently} of the electrons, and therefore models the
observed behaviour of neutrino oscillations between the three flavours. 

Thus the proposed splitting of the `generation' symmetry into two separate triplet symmetries is supported by experiment, at least to a certain extent. 
Also, there is a crucial
difference between the electron generation symmetry, in which the mass changes very obviously by several orders of magnitude between generations, and the neutrino generation
symmetry, in which no mass can be detected directly at all, but is instead inferred from the details of the theoretical model.
It is reasonable to suppose, therefore, that the electron generation symmetry acts to change the mass of the quarks, while the neutrino generation symmetry has
no detectable effect on quark mass.

 In other words, the neutrino generation symmetry effectively acts on the colours of the quarks rather than their generation.
But this implies that the neutrino generation should rather be interpreted as a colour, and that neutrino oscillation should be interpreted 
as a change in colour (as was implicit at the end of Section~\ref{particleids}).
This radical proposal has a number of interesting consequences. First, it allows neutrinos to change generation/colour by an interaction (with a massless gluon), 
rather than spontaneously,
and therefore it allows this to happen even if the neutrinos have exactly zero (rest) mass. 
It therefore does not require a hypothesis of non-zero neutrino mass.
Second, it allows colours to escape from hadrons, attached to neutrinos, and therefore removes the need to explain
`colour confinement'. Third, it allows free massless gluons to exist, which are then available as mediators for a quantum theory of gravity.
Indeed, it is possible that it even allows bound states of two spin $1$ gluons (with different colour charges) to form spin $2$ gravitons.

In order to build a mathematical model and make predictions to test such speculations, it is therefore necessary to introduce spinors
into the model, and Dirac matrices to act on the spinors, and to consider ways in which the three generations of neutrinos can be represented
mathematically within the spinor representations.

\section{Dirac matrices}
\subsection{The Dirac algebra}
In the first part of this paper, Sections~\ref{intro} to \ref{gauges}, we have been concerned only with the gauge group of the Standard Model,
and deriving its properties from the properties of a particular complex reflection group of order $648$. We have made no mention of spin or the Dirac matrices \cite{Dirac},
which are usually treated more-or-less independently of the gauge groups. However, the formalism of electro-weak mixing in the Standard Model
does involve a projection with $1-\gamma^5$, which effectively `mixes' part of the gauge group with the Dirac algebra describing spin. That is,
although $\gamma^5$ commutes with the relativistic spin group generated by even products of $\gamma^0,\gamma^1,\gamma^2$ and $\gamma^3$, 
as required by the Coleman--Mandula Theorem \cite{Coleman}, it does not commute with the odd products.

Moreover, if we are to explain the various equations given above and in \cite{finite}, that relate the mixing angles of the Standard Model to various particle masses,
then we must find some way to mix the Dirac equation (that defines the masses) with the gauge groups (that define the mixing angles). In other words, we must
mix the Dirac matrices with the Gell-Mann matrices and Pauli matrices into a single finite group that unifies all the finite symmetry operations.
If we are to integrate the Dirac algebra into the unification scheme proposed above, then we require a non-trivial action of the group $G_{27}\rtimes Q_8\rtimes Z_3$,
either on the Dirac algebra itself, or on a finite analogue of the Dirac algebra.
The obvious choice is to use the complex representation $\rep1b+\rep3b$ for the Dirac spinors, and tensor it with itself or its complex conjugate to get the
representation on the Dirac algebra. 

However, this does not actually work, for the simple reason that it does not contain a copy of $\rep2a$, which is required for the weak force to act on leptons.
It leads to a model of Pati--Salam \cite{PatiSalam} type, based on $SU(4)$, but does not unify $SU(4)$ with $SU(2)_L$. In fact we need a $4$-dimensional \emph{quaternionic}
representation, based on a complex representation something like $\rep2a+\rep3b+\overline{\rep3b}$.
It turns out to be easier to build the algebra first, and then get the group to act on it.
However, there is no action of the modified Gell-Mann matrices on the finite group generated by the Dirac matrices as they stand.
Therefore it is necessary to extend the Dirac matrices by adding at least one more generator.
But it is indeed possible, as we shall see shortly, to start with a finite group generated by 
six (rather than the standard five) Dirac matrices, 
and get the group $G_{27}\rtimes Q_8\rtimes Z_3$ to act on this group. 

In fact there is a unique possible action of the proposed finite analogue of the gauge group of the Standard Model
on an appropriate finite analogue of the Dirac matrices. This requires extending the Dirac matrices from complex to quaternionic,
in order to incorporate three generations of fermions, and extends the usual group of order $64$ to order $128$.
The resulting `unification' of Dirac, Gell-Mann and Pauli matrices into a single structure
forms a group
of order $82944$, consisting of $4\times 4$ quaternion matrices, that is a subgroup of an exceptional 
quaternionic reflection group \cite{quatref,BM}. This reflection group has a normal subgroup of order $2^7=128$, with quotient
group isomorphic to the rotation subgroup of the Weyl group of type $E_6$.
This fact may be the underlying reason why many people have seen an imprint of $E_6$ symmetry in the Standard Model,
although here we have used the finite groups to break the $E_6$ symmetry down to something resembling the symmetries we see in particle experiments.

If this group is useful for physics, then it certainly goes beyond the Standard Model of Particle Physics, because it allows
the Gell-Mann matrices and the Pauli matrices to act on the Dirac matrices, instead of commuting with them. In other words,
it allows particle physics to determine, or at least to influence, the shape of spacetime, and possibly thereby to include a quantum theory of gravity.
It is not at all obvious what such a theory should look like, but
if we develop the mathematical model first, then the physical mechanisms may become clearer later.

\subsection{A $4$-dimensional group}
\label{Diracmats}
The Dirac matrices are $4\times 4$ complex matrices that have a similar structure to that of the Pauli matrices and Gell-Mann matrices.
But in this case both Hermitian and anti-Hermitian matrices are required, and the group that is normally used is $SL(2,\CC)$, rather than
$SU(4)$ or $SL(4,\RR)$. The Dirac matrices are all unitary, and they generate a group of order $2^6=64$, that is a commuting product of $Z_4$
with two copies of $Q_8$. 
A more natural analogue of $Q_8$ and $G_{27}$, however, is a so-called \emph{extraspecial} group $E_{128}$ of order $128$ obtained by adjoining a complex conjugation
operator.

In principle, there is a choice as to whether the complex conjugation operator squares to $+1$ or $-1$, 
because there are two extraspecial groups of order $128$. One is the tensor product (also called central product) of three copies of the quaternion group $Q_8$,
while the other is the tensor product of three copies of the dihedral group $D_8$. In order to get an action of $G_{27}$ on this group,
it is necessary to have automorphisms of order $3$
of each factor, and therefore to take $E_{128}$ as the tensor product of three copies of $Q_8$.
In this case, the complex conjugation operator squares to $-1$.
Then the representation is quaternionic, which means it can also be written as
a group of $4\times 4$ quaternion matrices. 

The group $E_{128}$ requires six generators, such as 
\begin{align}\label{E128gens}
&\begin{pmatrix}i&0&0&0\cr 0&i&0&0\cr 0&0&i&0\cr 0&0&0&i\end{pmatrix}, \quad
\gamma^3=\begin{pmatrix}0&0&1&0\cr 0&0&0&-1\cr -1&0&0&0\cr 0&1&0&0\end{pmatrix},\quad
i\gamma^1\gamma^2\gamma^3 =\begin{pmatrix}0&0&-1&0\cr 0&0&0&-1\cr 1&0&0&0\cr 0&1&0&0\end{pmatrix},\cr
&\begin{pmatrix}j&0&0&0\cr 0&j&0&0\cr 0&0&j&0\cr 0&0&0&j\end{pmatrix}, \quad
\gamma^1=\begin{pmatrix}0&0&0&1\cr 0&0&1&0\cr 0&-1&0&0\cr -1&0&0&0\end{pmatrix},\quad
i\gamma^0\gamma^2=\begin{pmatrix}0&0&0&1\cr 0&0&-1&0\cr 0&1&0&0\cr -1&0&0&0\end{pmatrix}.
\end{align}
These generators are arranged in such a way that the three columns generate three mutually commuting copies of the quaternion group $Q_8$,
and the rows generate abelian subgroups isomorphic to $Z_4\times Z_2\times Z_2$.

The names above are given according to the standard Bjorken-Drell convention \cite{BjorkenDrell}, but many other arrangements are possible. For example, it is possible to choose
$\gamma^1\gamma^2$, $\gamma^2\gamma^3$ and $\gamma^1\gamma^3$ as elements of one of the $Q_8$ factors, and $i\gamma_0$, $i\gamma_5$ and
$\gamma_0\gamma_5$ as elements of another $Q_8$ factor. This choice gives a mathematical separation of the properties that depend on a
direction in space (such as spin and momentum) from those that do not (such as energy, mass and weak isospin), which may be useful for physical purposes.

It is well-known that the outer automorphism group of this group of order $2^7$ is isomorphic to the Weyl group of type $E_6$, so that it is
straightforward to calculate any necessary properties.
In fact there are $120$ subgroups isomorphic to $Q_8$, and they are all equivalent under the automorphism group, so we can take any one of them to represent the
non-relativistic spin group as above. We then need to choose $\gamma^0$ and $\gamma^5$ from the nine elements of order $2$ (up to sign) that commute with this
copy of $Q_8$. The automorphism group acts transitively on these $9$ elements, and the stabiliser of one of them acts transitively on the four that anti-commute with it.
Hence all choices of Dirac matrices are mathematically equivalent, as long as we still have the full $E_6$ symmetry.

\subsection{A quaternionic reflection group}

The action of $G_{27}$ on $E_{128}$ is more difficult to calculate, but much of the work has been done for us in \cite{quatref,BM},
though with some errors in \cite{BM} we need to watch out for.
In order to avoid overloading the paper with detailed calculations, I shall state without proof numerous facts about this action that are well-known to
group-theorists. 
It is probably easiest to start from the quaternionic reflection group called $S_3$ in \cite{quatref},
which happens to be an extension of $E_{128}$ by the rotation subgroup of the Weyl group of type $E_6$. Quaternionic reflections are given by 
the same formula (\ref{reflform}) as for complex reflections, with appropriate choices for the quaternionic scalar $\lambda$.

In order to ensure that reflections are linear maps
we should check the order of multiplication, since the quaternions are non-commutative. 
(Note that the conventions in \cite{quatref,BM} are opposite to each other, and we follow the latter.) 
The particular order given in (\ref{reflform}) ensures that replacing $x$ by a scalar multiple $\mu x$ gives the same scalar multiple on the right hand side.
Moreover, replacing $r$ by a scalar multiple $\mu r$ replaces $x.r$ by $(x.r)\bar\mu$, so changes the scalar $1-\lambda$ to $1-\mu^{-1}\lambda\mu$.
Hence reflections of order $2$, with $\lambda=-1$, are invariant under this operation, but other reflections are not.

The group $G_{27}$ can be generated by the products of pairs of reflections of order $2$ in the mirrors defined by
\begin{align}
(2,0,0,0), (\pm1,1,1,1), (\pm1,i,j,k), (\pm1,j,k,i), (\pm1,k,i,j).
\end{align}
Indeed, we can take the products of the first reflection with the other $8$ to represent the $8$ Gell-Mann matrices.
Generators for the full group of order $648$ 
can then be taken as 
\begin{align}
\frac12\begin{pmatrix}-1&1&1&1\cr -1&1&-1&-1\cr -1&-1&1&-1\cr -1&-1&-1&1 \end{pmatrix},\quad
\begin{pmatrix}1&0&0&0\cr 0&i&0&0\cr 0&0&j&0\cr 0&0&0&k\end{pmatrix}, \quad
\begin{pmatrix} 1&0&0&0\cr 0&0&1&0\cr 0&0&0&1\cr 0&1&0&0\end{pmatrix}
\end{align}
The first matrix above represents one of the Gell-Mann matrices in $G_{27}$, the second one of the Pauli matrices in $Q_8$, and the last 
is a generator for $Z_3$.

But now we must decide on a labelling of the Dirac matrices with $\gamma^\mu$, since the action of the Gell-Mann matrices breaks the
symmetry of the $120$ copies of $Q_8$. The full normalizer of $G_{27}$ is $G_{27}\rtimes Q_8\rtimes Z_3$, which breaks $120$ into $12+108$,
and therefore breaks the symmetry between the standard (Bjorken--Drell or chiral) labelling, which belong to the set of $108$, and the more 
symmetrical version suggested above, which is one of the $12$. 
It would seem that these two labellings imply physically different
interactions between the colours/generations and spin/spacetime/gravity, and therefore at most one of them can provide a correct model of physics.

However, the real issue is the correct interpretation of all the groups. For example, do we want a model that is invariant under an $SO(3)$ local rotation group,
or do we want a model that couples to gravity, and breaks the symmetry of space by reference to the direction of the gravitational field?
If the former, then we must use the new symmetrical labelling of Dirac matrices. If the latter, then we can stick to the traditional labellings. 
We may therefore need to keep both possibilities in mind, and potentially use both in different contexts.

\subsection{Clifford algebras}
The full $4\times 4$ quaternion matrix algebra can also be interpreted as a Clifford algebra for a $6$-dimensional real space with signature $(6,0)$,
$(5,1)$, $(2,4)$ or $(1,5)$. Many models of this type have been proposed \cite{PatiSalam,Stoica,twistors,twistorlectures}, 
including the Pati--Salam model based on $SU(4)=Spin(6)$,
and the twistor theory of Penrose,
based on the Clifford algebra $Cl(2,4)$ and the corresponding spin group isomorphic to $SU(2,2)$. Twistors can equally well be given the opposite
signature, so that the matrix algebra underlying $Cl(4,2)$ is the algebra of $8\times 8$ real matrices. Models of this type have also been proposed \cite{WoitSO24},
as have models based on $Cl(3,3)$, and it is not at all clear which of these Clifford algebras 
provides the closest fit to the Standard Model,
and/or to observed physics. 

The Dirac matrices above now act on rows (mathematicians' convention) or columns (physicists' convention)
of four quaternions, as a generalisation to three generations of the standard Dirac spinors,
that are columns of four complex numbers. The projections with $1\pm \gamma^5$ onto left-handed and right-handed (Weyl) spinors
now project onto quaternionic $2$-spinors, that contain some kind of generation/colour label $i$, $j$, or $k$ in addition to
the standard interpretation as spin.
At the same time, we see that three generations only require twice as many spinors as one generation,
as has been pointed out in \cite{octions,chirality}. 
It is usually assumed (without proof) that three generations require three times as many spinors as one
\cite{DG}.

\subsection{Twistors}
The Dirac generators $\gamma^0$, $\gamma^1$, $\gamma^2$, $\gamma^3$ and $\gamma^5$ generate the Clifford algebra $Cl(2,3)$.
Hence in order to obtain $Cl(2,4)$ it is necessary and sufficient
to add an extra generator $\gamma^6$ that anti-commutes with all of the five generators, and squares to $-1$.
In both the Bjorken--Drell and chiral conventions 
we have 
\begin{align}
\gamma^6 = \begin{pmatrix}0&0&0&j\cr 0&0&-j&0\cr 0&-j&0&0\cr j&0&0&0\end{pmatrix}.
\end{align}
However, this choice breaks the symmetry of the three generations $i$, $j$ and $k$, 
as well as the $SO(3)$ symmetry of space, since $\gamma^2$ has
a generation label $i$, but $\gamma^1$ and $\gamma^3$ have not. The former picks out the first generation of `normal' matter, while the latter
conceivably picks out the direction of the gravitational field. 

The twistors themselves are much 
the same as the quaternionic Dirac spinors, and can be represented as quaternionic $4$-vectors, either as rows or as columns
according to convention. However, the quaternionic structure is not used in most approaches to twistors, which are considered to lie
in spinor representations of $Spin(2,4)\cong SU(2,2)$, and are therefore complex $4$-vectors. Each quaternionic Dirac spinor
therefore consists of a chiral pair of Penrose twistors.

For convenience I list the generators in the two conventions. First I give the standard chiral convention, in which
the spin symmetries are broken:
\begin{align}\label{brokenlabels}
&\gamma^1=\begin{pmatrix}0&0&0&1\cr 0&0&1&0\cr 0&-1&0&0\cr -1&0&0&0\end{pmatrix},
\gamma^2=\begin{pmatrix}0&0&0&-i\cr 0&0&i&0\cr 0&i&0&0\cr -i&0&0&0\end{pmatrix},
\gamma^3=\begin{pmatrix}0&0&1&0\cr 0&0&0&-1\cr -1&0&0&0\cr 0&1&0&0\end{pmatrix},\cr
&\gamma^0=\begin{pmatrix}0&0&1&0\cr 0&0&0&1\cr 1&0&0&0\cr 0&1&0&0\end{pmatrix},
\gamma^5=\begin{pmatrix}1&0&0&0\cr 0&1&0&0\cr 0&0&-1&0\cr 0&0&0&-1\end{pmatrix},
\gamma^6 = \begin{pmatrix}0&0&0&j\cr 0&0&-j&0\cr 0&-j&0&0\cr j&0&0&0\end{pmatrix}.
\end{align}
For comparison, here is an example of the alternative `symmetric' convention
in which the spin symmetries are unbroken:
\begin{align}\label{unbrokenlabels}
&\gamma^1=\begin{pmatrix}0&i&0&0\cr i&0&0&0\cr 0&0&0&i\cr 0&0&i&0\end{pmatrix},
\gamma^2 =\begin{pmatrix}0&0&0&i\cr 0&0&-i&0\cr 0&-i&0&0\cr i&0&0&0\end{pmatrix},
\gamma^3=\begin{pmatrix}i&0&0&0\cr 0&-i&0&0\cr 0&0&i&0\cr 0&0&0&-i\end{pmatrix},\cr
&\gamma^0=\begin{pmatrix}0&0&0&i\cr 0&0&-i&0\cr 0&i&0&0\cr -i&0&0&0\end{pmatrix},
\gamma^5 = \begin{pmatrix}0&i&0&0\cr -i&0&0&0\cr 0&0&0&-i\cr 0&0&i&0\end{pmatrix},
\gamma^6 = \begin{pmatrix}j&0&0&0\cr 0&j&0&0\cr 0&0&j&0\cr 0&0&0&j\end{pmatrix}.
\end{align}
It remains to be seen which of the two versions provides a better description of observed physics in a gravitational field.
Let us call these two labellings `broken' and `unbroken' respectively, for ease of reference.

\subsection{Action of Gell-Mann matrices on Dirac matrices}
In both cases, we can translate to a representation of the Gell-Mann matrices by mapping each copy of $Q_8$ to $\FF_4$, via the map that takes
$\pm1$, $\pm i$, $\pm j$ and $\pm k$ to $0$, $1$, $v$ and $w$ respectively. 
Thus the tensor product of three copies of (multiplicative) $Q_8$ translates to the direct sum of three copies of (additive) $\FF_4$.
There is some choice as to which elements of each $Q_8$ are labelled $1$ and $v$ and $w$,
so only the general shape is really important.
In the standard chiral convention we have
\begin{align}
\begin{array}{c|c}
\gamma^1\sim (0,0,w) & \gamma^0 \sim (0,v,1)\cr
\gamma^2\sim (1,1,1) & \gamma^5\sim (0,w,1)\cr
\gamma^3\sim (0,0,v) & \gamma^6\sim (v,1,1)
\end{array}
\end{align}
while in the proposed symmetric convention we have
\begin{align}
\begin{array}{c|c}
\gamma^1\sim (1,1,1) & \gamma^0\sim (1,v,0)\cr
\gamma^2\sim (1,1,v) & \gamma^5\sim (1,w,0)\cr
\gamma^3\sim (1,1,w) & \gamma^6\sim (v,0,0)
\end{array}
\end{align}

Multiplication of Dirac matrices corresponds to addition of vectors over $\FF_4$, so that the Lorentz group in the chiral convention has 
generators
\begin{align}
&(1,1,v), (1,1,w), (0,0,1),\cr
&(0,v,v), (1,w,0), (0,v,w), 
\end{align}
in which the first row contains the rotations and the second row the boosts.
Again, the breaking of spatial symmetry is clear from the notation alone.
In the symmetric convention, we have generators
\begin{align}
 &(0,0,1), (0,0,v), (0,0,w),\cr
&(0,w,1), (0,w,v), (0,w,w). 
\end{align}
It is clear that, in order to have a model of the strong force that is invariant under rotations in space, we must switch to the
symmetric convention for labelling Dirac matrices. But that implies a wholesale relabelling of the Dirac spinors as well, which is the subject of the next section.

Next we consider compatibility with the weak interaction, defined by the actions of the Pauli matrices in (\ref{ternaryPauli}). Here the fixed space of the Pauli matrices
is spanned by $(0,1,-1)$ in the complex $3$-space, or $(0,1,1)$ in the finite $3$-space. If we require the weak interaction to be invariant under space rotations, then 
this one-space must be fixed by the space rotation group, which implies that we have to re-order the coordinates of the complex $3$-space,
say to the following:
\begin{align}
\begin{array}{c|c}
\gamma^1\sim (1,1,1) & \gamma^0\sim (0,v,1)\cr
\gamma^2\sim (v,1,1) & \gamma^5\sim (0,w,1)\cr
\gamma^3\sim (w,1,1) & \gamma^6\sim (0,0,v)
\end{array}
\end{align}
On the other hand, the known correlation between the weak interaction and the direction of spin, as demonstrated by the Wu experiment \cite{Wu},
may be better modelled in the original coordinate system. Essentially, the question comes down to what type of model we want, and in particular whether we want
a model in which gravity is independent of the other forces, or a model in which gravity is unified with, and therefore couples to, the other forces. 

\subsection{$C$, $P$ and $T$ symmeties}
In any case, the three copies of $Q_8$ into which $E_{128}$ factorises correspond to the three coordinates of the complex $3$-dimensional representation, and in the
symmetric (unbroken) convention are as follows:
\begin{align}
&\langle \gamma^1\gamma^2, \gamma^2\gamma^3, \gamma^3\gamma^1\rangle,\cr
&\langle i\gamma^0, i\gamma^5,\gamma^0\gamma^5\rangle,\cr
&\langle i,j,k\rangle.
\end{align}

In each case, negation of all 
the elements in $Q_8$ is an anti-automorphism that reverses the order of multiplication. It corresponds to
complex conjugation, since it swaps $v=(-1+i+j+k)/2$ with $w=(-1-i-j-k)/2$.
In the first case, this is the parity symmetry $P$, that is a discrete version of negating the three coordinates of space.
In the second case, it negates the `time' coordinate $i\gamma_0$, and is therefore a version of the time-reversal symmetry $T$.
In the third case, it negates $i$, and therefore acts as a charge-conjugation symmetry $C$, in this case acting on all three generations,
labelled by $i$, $j$ and $k$. No combination of $P$, $T$ and $C$ is a symmetry of the finite group, but the combination $CPT$ is a
symmetry of the representations, since it swaps $\rep3b$ with $\overline{\rep3b}$. 

Hence there is a sense in which the model
has $CPT$ symmetry,
but
of course, $CPT$ is not a \emph{physical} symmetry. 
The symmetry of the model merely suggests that $CPT$-violation cannot be detected experimentally, at least if gravity is not taken into account.
On the other hand, violations of the individual symmetries $P\equiv CT$, $T\equiv CP$ and $C\equiv PT$ can all be detected experimentally. The first was
verified by the Wu experiment \cite{Wu} in 1957, that detected a correlation between the direction of spin of a cobalt 60 nucleus and the direction of
the electron ejected in beta decay. 
The second was confirmed by the Cronin--Fitch experiment \cite{CP} in 1964, that detected a mixing between different neutral kaon decay modes, such that a kaon beam that 
was supposedly a pure beam that should decay into an odd number of pions actually exhibited a small proportion of even-pion decay. 
The third is an asymmetry between matter and anti-matter, and is 
generally considered to be obvious (from an experimental point of view).

From a theoretical point of view, on the other hand, there is no generally accepted explanation for any of these three types of fundamental symmetry-breaking.
Our proposal for mixing of the gauge groups with the Dirac matrices provides a possible mechanism for explaining one or more of them in terms of quantum gravity.

\subsection{The Lorentz group}
\label{Lorentzgp}
The first two copies of $Q_8$ generate a copy of the Lie group $SO(4)$, consisting of two commuting copies of $SU(2)$. The $9$ products of one generator
from each copy of $Q_8$ are
\begin{align}
&\gamma^0\gamma^1, \gamma^0\gamma^2, \gamma^0\gamma^3,\cr
&\gamma^5\gamma^1,\gamma^5\gamma^2,\gamma^5\gamma^3,\cr
&i\gamma^1,i\gamma^2,i\gamma^3.
\end{align}
The first three extend $SU(2)$ to the standard copy of the Lorentz group $SL(2,\CC)$,
while each of the other rows extends $SU(2)$ to a different copy of $SL(2,\CC)$.
The three rows together extend $SO(4)$ to $SL(4,\RR)$.
The question then is how best to interpret these groups. Since $SL(2,\CC)$ is normally interpreted as the double cover of the
Lorentz group $SO(3,1)$, it is natural to interpret $SL(4,\RR)$ as the double cover of $SO(3,3)$, and therefore postulate a six-dimensional
`spacetime', as in \cite{Chester}. However, in the dual formulation of the Lorentz group, acting on momentum-energy, there is a more obvious
interpretation as three dimensions of momentum and three of angular momentum. This allows us to keep a four-dimensional spacetime,
as required by physical considerations, and, with a scalar mass, represent momentum and angular momentum in the antisymmetric square
of the spacetime representation.

At the same time, this interpretation breaks the Hamiltonian duality between spacetime and $4$-momentum, by insisting that energy must always have
a direction in space. It thereby breaks the standard identification of $SL(2,\CC)$ with $SO(3,1)$, which now become disjoint subgroups of $SL(4,\RR)$.
Indeed, the identification of these two groups was originally only an analogy, using the fact that the two groups have the same Lie algebra,
and is still to this day presented in textbooks as an analogy,
not as a mathematically rigorous deduction. Downgrading this identification to its original status as an analogy affects only interpretations,
and not the calculations that are essential to quantum mechanics and particle physics. The main benefit of such downgrading is the upgrading
of the Standard Model of Particle Physics from a Lorentz-covariant theory to a generally covariant theory. 
This 
would remove the biggest single obstacle to unification of particle physics with gravity.

\section{The structure of the combined group}
\label{structure}

\subsection{Gell-Mann matrices}
The full set of $8$ Gell-Mann matrices can be taken as the following four matrices and their conjugate transposes:
\begin{align}\label{GellMann4Q}
\frac12\begin{pmatrix}-1&1&1&1\cr -1&1&-1&-1\cr -1&-1&1&-1\cr -1&-1&-1&1 \end{pmatrix},\quad
\frac12\begin{pmatrix}-1&i&j&k\cr i&1&k&-j\cr j&-k&1&i\cr k&j&-i&1 \end{pmatrix},\cr
\frac12\begin{pmatrix}-1&j&k&i\cr j&1&i&-k\cr k&-i&1&j\cr i&k&-j&1 \end{pmatrix},\quad
\frac12\begin{pmatrix}-1&k&i&j\cr k&1&j&-i\cr i&-j&1&k\cr j&i&-k&1 \end{pmatrix}.
\end{align}
The two real matrices here correspond to the two colourless gluons, and the other six to the coloured gluons. Writing $v=(-1+i+j+k)/2$,
the group of order $27$ generated by these matrices
contains a centre of order $3$, generated by
\begin{align}\label{G27centre}
\begin{pmatrix}v&0&0&0\cr 0&0&v&0\cr 0&0&0&v\cr 0&v&0&0\end{pmatrix}.
\end{align}

This matrix acts on the Dirac matrices to permute the three generations defined by $i$, $j$ and $k$. 
It is important to note that replacing $v$ by $w$ in (\ref{G27centre}) can only be done if the coordinate permutation is also inverted.
The chirality of the weak interaction ultimately
arises from this important fact, as we shall see later. 
It is also worth noting that this matrix fixes the quaternionic one-space spanned by $(0,1,v,w)$, on which all the Gell-Mann matrices act trivially.

\subsection{Pauli matrices} 
The (anti-Hermitian) Pauli matrices generating $Q_8$ are represented by the diagonal matrices
\begin{align}\label{Pauli4Q}
\begin{pmatrix}1&0&0&0\cr 0&i&0&0\cr 0&0&j&0\cr 0&0&0&k\end{pmatrix}, \quad
\begin{pmatrix}1&0&0&0\cr 0&j&0&0\cr 0&0&k&0\cr 0&0&0&i\end{pmatrix}, \quad
\begin{pmatrix}1&0&0&0\cr 0&k&0&0\cr 0&0&i&0\cr 0&0&0&j\end{pmatrix}, 
\end{align}
acting on the spinors. These matrices commute with (\ref{G27centre}), but are permuted by both the `scalar' $v$ and the coordinate permutation separately.
The scalar acts by permuting the three generations, but leaving the spin direction alone, while the permutation acts on the spin direction but not on the generation.
The fact that the matrix ({\ref{G27centre}) couples the cyclic ordering of the three generations to the cyclic ordering of the three directions of spin is what is
known as the \emph{chirality} of the weak interaction. Conventionally this chirality is \emph{left-handed}, so that the copy of $SU(2)$ generated by the
Pauli matrices above is denoted $SU(2)_L$.

The way chirality is implemented in the Standard Model is by restricting to the third component of weak isospin, represented say by the first matrix in (\ref{Pauli4Q}),
and realising $j$ as complex conjugation, so that the bottom half of the Dirac spinor behaves as a Weyl spinor that is the complex conjugate of the top Weyl spinor.
Thus the proposed model is a straightforward generalisation of the Standard Model at this point, and incorporates all three components of weak isospin,
not just one.

An alternative interpretation of the three components of weak isospin is to take the standard `third component' to be represented by the sum of the three matrices
(\ref{Pauli4Q}), so that it corresponds to the complex number $i+j+k$ in all three directions in space, and thereby allows the weak force to be generally covariant.
This interpretation is consistent with the analysis of charge, weak hypercharge and weak isospin in Section~\ref{theBTG}.

We can now extend from $Q_8$ to the binary tetrahedral group by adjoining any one of the three matrices
\begin{align}\label{U1actions}
\begin{pmatrix}1&0&0&0\cr 0&0&1&0\cr 0&0&0&1\cr 0&1&0&0\end{pmatrix},\quad
\begin{pmatrix}v&0&0&0\cr 0&v&0&0\cr 0&0&v&0\cr 0&0&0&v\end{pmatrix},\quad
\begin{pmatrix}v&0&0&0\cr 0&0&0&v\cr 0&v&0&0\cr 0&0&v&0\end{pmatrix}.
\end{align}
As noted above, 
these three copies of the binary tetrahedral group are not equivalent, so that the physical interpretations of these three different $Z_3$ symmetry groups
must be different. They act on various triplet symmetries including generations, colours, and directions of spin, in various combinations yet to be determined.

\subsection{Classification of spinors}
We have identified four distinct copies of $U(1)$, generated by the following four matrices in $U(3)$:
\begin{align}
\begin{pmatrix}v&0&0\cr 0&v&0\cr 0&0&v\end{pmatrix}, \quad
\begin{pmatrix}v&0&0\cr 0&1&0\cr 0&0&1\end{pmatrix}, \quad
\begin{pmatrix}1&0&0\cr 0&v&0\cr 0&0&v\end{pmatrix}, \quad
\begin{pmatrix}v&0&0\cr 0&w&0\cr 0&0&w\end{pmatrix}.
\end{align} 
The first is a scalar, the second is a complex reflection, and the other two are products of a scalar and a reflection. They all have different physical interpretations,
and they all act differently on the spinors. Here we use mathematicians' convention, writing spinors as rows, with scalars on the left, and matrix action on the right.
To switch to physicists' convention, transpose the matrices, write spinors as columns, scalars on the right, and matrix action on the left.
The scalar acts on spinors as in (\ref{G27centre}), and acts trivially on the spinors that are (left) scalar multiples of $(0,1,v,w)$. 

The other three copies of $U(1)$ act via the matrices in (\ref{U1actions}), in which the first matrix represents the complex reflection.
The last matrix acts trivially on just one spinor, consisting of the scalar multiples of $(0,1,w,v)$, while the second one acts non-trivially on all spinors.
The first matrix acts trivially on the spinors generated by $(1,0,0,0)$ and $(0,1,1,1)$, which can therefore be mixed in arbitrary (quaternionic) proportions
to create something like the left- and right-handed Weyl spinors of the Standard Model. For example, one could take the spinors generated by
$(i+j+k,1,1,1)$ and $(-i-j-k,1,1,1)$, and use the first coordinate as an imaginary time coordinate, and the other three as real space coordinates.
This interpretation is consistent with the discussion in Section~\ref{Lorentzgp},
in which the four quaternionic coordinates are acted on by $SL(4,\RR)$ as a group of spacetime symmetries derived from General Relativity.

Alternatively, one can keep time and space separate, so that the Pauli matrices act trivially on $(1,0,0,0)$, and non-trivially on all the other spinors.
This would allow us to interpret $(1,0,0,0)$ as a right-handed spinor, and the rest as left-handed. Moreover, $(0,1,v,w)$ is acted on as a quaternionic one-space,
which could be interpreted as a left-handed spinor for leptons,     
leaving $(0,1,1,1)$ and $(0,1,w,v)$ as left-handed spinors for quarks. The fact that there are only two independent spinors for quarks rather than three
is then equivalent to colour confinement.

\subsection{Action of Gell-Mann matrices} 
We now consider the representation of the Gell-Mann matrices (\ref{GellMann4Q}). The first pair, that is the real matrices, representing colourless gluons,
commute with the `flavour' $Q_8$, consisting of `scalar' matrices. Since they have order $3$, they cannot swap the other two factors, so they fix
them both. The following calculation shows that they cycle $\gamma^3$, $\gamma^1$ and $\gamma^1\gamma^3$ in the standard convention:
\begin{align}
\frac12\begin{pmatrix}-1&-1&-1&-1\cr 1&1&-1&-1\cr 1&-1&1&-1\cr 1&-1&-1&1 \end{pmatrix}
\begin{pmatrix}0&0&1&0\cr 0&0&0&-1\cr -1&0&0&0\cr 0&1&0&0\end{pmatrix}
&\frac12\begin{pmatrix}-1&1&1&1\cr -1&1&-1&-1\cr -1&-1&1&-1\cr -1&-1&-1&1 \end{pmatrix}\cr
&= \begin{pmatrix}0&0&0&1\cr 0&0&1&0\cr 0&-1&0&0\cr -1&0&0&0\end{pmatrix}
\end{align}
It is hard to interpret this in a physically reasonable way, except to say that on the scale of the strong force, the structure of spacetime
completely breaks down. With the proposed symmetric labelling, on the other hand, this calculation implies a cycling of $\gamma^1\gamma^2$ with
$\gamma^2\gamma^3$ and $\gamma^3\gamma^1$, 
that is, a change in the direction of spin.

The action on the third copy of $Q_8$ is similar:
\begin{align}
\frac12\begin{pmatrix}-1&-1&-1&-1\cr 1&1&-1&-1\cr 1&-1&1&-1\cr 1&-1&-1&1 \end{pmatrix}
\begin{pmatrix}0&0&-1&0\cr 0&0&0&-1\cr 1&0&0&0\cr 0&1&0&0\end{pmatrix}
&\frac12\begin{pmatrix}-1&1&1&1\cr -1&1&-1&-1\cr -1&-1&1&-1\cr -1&-1&-1&1 \end{pmatrix}\cr
&=
\begin{pmatrix}0&-1&0&0\cr 1&0&0&0\cr 0&0&0&1\cr 0&0&-1&0\end{pmatrix}
\end{align}
Again, this seems to have no reasonable physical meaning in the standard labelling, but in the proposed symmetric labelling
we see a conversion from $\gamma^0\gamma^5$ to $i\gamma^5$ (and also to $i\gamma^0$), which indicates a mixing of some kind with the
weak interaction. 

\subsection{Colour}
The action of the coloured gluons is a good deal more complicated, since they act by permuting the three copies of $Q_8$. Thus they map
spin symmetries to other types of internal symmetries, 
which seems to imply that the entire structure of spacetime breaks down once we look at the internal structure of a baryon.
In other words, at distances less than about $10^{-15}$m, it no longer makes sense to talk about `space' at all. Dividing by the speed of light to
convert distance to time, we obtain a cutoff of around $3\times 10^{-24}$s below which it does not make sense to talk about `time'. This is 
roughly half the lifetime of the most unstable mesons and baryons in the Standard Model.

It is straightforward to 
compute the action of the coloured gluons on the Dirac matrices, and verify that the action is compatible with (\ref{GellMannunitary}).
Hence the latter notation provides a simpler way to see the action, if we interpret the three complex coordinates as representing the three
$Q_8$ factors. In this way the diagonal matrices represent the colourless gluons, and the coloured gluons permute the three factors.

\section{Intimations of quantum gravity}
\subsection{Choosing a gauge group}
In the Standard Model, gauge groups are obtained by a process of exponentiating Hermitian (or anti-Hermitian, depending on convention) matrices.
This process still works in the case of the weak gauge group $SU(2)$, but in the case of the strong gauge group $SU(3)$ our matrices are no longer
Hermitian or anti-Hermitian, so the process needs modification. What is really going on at a fundamental level is that the gauge group is
the smallest Lie group that contains the finite group, in its fundamental representation. With this definition the gauge group is automatically compact, and unitary, but is not
necessarily the full unitary group. If we take the finite group to be either $Q_8$ or $Q_8\rtimes Z_3$ in the case of the weak force, and either $G_{27}$ or
$G_{27}\rtimes Q_8$ in the case of the strong force,
then the gauge groups defined in this way are $SU(2)$ and $SU(3)$ respectively, in agreement with the Standard Model.

If we adopt the same principle in the case of the group $E_{128}$, things get a bit more complicated, since $E_{128}$ itself lies in a tensor product of three
 commuting copies of $SU(2)$, but $E_{128}\rtimes G_{27}$ does not. There are therefore various different groups that might be considered to be the 
 appropriate `gauge group'  of quantum gravity, depending on which properties of the group we consider to be most desirable. For example, in the cases of the weak and
 strong forces, the dimension of the Lie group is one less than the order of the derived quotient of the finite group:
\begin{align}
Q_8/Z_2 & = Z_2 \times Z_2,\cr
G_{27}/Z_3 & = Z_3 \times Z_3.
\end{align}
This principle would give us a dimension of $63$ in the case
\begin{align}
E_{128}/Z_2 = Z_2 \times Z_2 \times Z_2 \times Z_2 \times Z_2 \times Z_2,
\end{align}
and therefore a gauge group $SL(4,\HH)$, which is not compact.

We might consider compactness to be an essential property, in which case we might restrict to the compact subgroup $Sp(4)$, of dimension $36$.
Or we might consider non-compactness to be an essential property, in order to relate to the group $GL(4,\RR)$ that defines general covariance
in the theory of relativity. Or we might consider Lorentz-invariance of the gauge group to be an essential property, in which case we are restricted to either $SU(2)$, if we
demand compactness and/or general invariance, or $SL(2,\CC)$, if we allow (or demand) non-compactness. The finite group model suggests
that all of these are viable approaches, and that the choice between them is simply a choice of what type of model we want to build.

For a totally quantised `theory of everything', the compact group $Sp(4)$ of dimension $36$ is most appropriate. It contains the $12$-dimensional
gauge group $U(1)\times SU(2)\times SU(3)$ of the Standard Model, plus $24$ further dimensions that contain what are, from the Standard Model point of view,
effectively arbitrary parameters. For a theory that combines the Standard Model with General Relativity, while making the minimal
possible changes to both, consistent with eliminating the known contradictions between them, we should reduce from the non-compact group
$SL(4,\HH)$ to its subgroup $SU(2)\otimes SL(4,\RR)$, so that \emph{locally} the theory of particle physics decouples from gravity as described by
General Relativity. In this case, $SU(2)$ acts as a gauge group for quantum gravity in the traditional Yang--Mills sense, while 
the extra symmetries of $SL(4,\HH)$ create a mixing between particle physics and gravity,
that can only be detected by \emph{non-local} experiments, .

In other words, we should utilise the splitting of $E_{128}$ as a commuting product of three copies of $Q_8$, and only exponentiate the elements that lie in one (or perhaps two)
of the factors. Indeed, two copies of $Q_8$ give a $15$-dimensional gauge group $SL(4,\RR)$, containing the Lorentz group $SL(2,\CC)$, and with compact part
$SO(4)$. The latter splits into two commuting copies of $SU(2)$, that are sometimes called left-handed and right-handed \cite{WoitSO24}. It is plausible that these
groups can be used to reproduce General Relativity, in which case the third copy of $Q_8$ provides a third gauge group $SU(2)$ that is not only Lorentz-invariant,
but generally invariant. This is the copy generated by the quaternionic scalars $i$ and $j$.

This copy of $SU(2)$ is not in the Standard Model, and is not in General Relativity. It therefore represents the (only) addition that is required in order to unify
particle physics and gravity into a single theory. As a gauge group, it gauges three quantities that link particle physics to gravity. Since the generation symmetry
acts on these three quantities, it is reasonable to suppose that they can be taken as the masses of the three generations of electron, relative to a standard unit of mass,
say the proton or the neutron. What it must do, therefore, is to reconcile (or calibrate) the inertial masses that are measured in particle physics
with the (active) gravitational masses that are measured with Cavendish-type experiments \cite{Cavendish,Gillies,newG,QLi}.

A similar mathematical idea is proposed in \cite{Chester}, where 
a group $Spin(3,3)\cong SL(4,\RR)$ is used to incorporate three generations of fermions by extending the concept of `time'
to three dimensions. The difference really is in the interpretation, and in particular the question of which mathematical concept do we label with the physical name
of `spacetime'. For compatibility with General Relativity, we need to label the $4$-dimensional representation of $SL(4,\RR)$ as spacetime, and not the $6$-dimensional
representation. Of course, this labelling is inconsistent with the conventional identification of $SL(2,\CC)$ as the Lorentz group $SO(3,1)$, or more precisely
the connected component of the identity $SO(3,1)^o$, acting on spacetime. But it is required for general covariance of any model of particle physics,
for the simple reason that $SL(4,\RR)$ does not have a double cover that extends the double cover $SL(2,\CC)$ of $SO(3,1)$.
Indeed, the real part of the quaternionic $4$-space is the tensor product of spinor representations of two copies of $SU(2)$, so it is most naturally
interpreted as a vector, and not a spinor.

\subsection{Gravi-weak mixing}
If we do take the copy of $SU(2)$ generated by the scalar matrices $i$, $j$ and $k$ as the gauge group $SU(2)_G$ for quantum gravity, as suggested above,
then there is some interesting mixing with the weak copy $SU(2)_W$, generated by the diagonal matrices with entries $(1,i,j,k)$, $(1,j,k,i)$ and $(1,k,i,j)$.
These two groups do not commute, and the finite group generated by the two copies of $Q_8$ in fact has order $2^8$. The commutator subgroup consists of
all sign changes on the four quaternion coordinates, and modulo signs we have four scalars plus $12$ even permutations of $(1,i,j,k)$.

This group therefore distinguishes four types of `spin', corresponding to the four quaternionic coordinates. In a gravitational context, these four coordinates are
labelled by the four coordinates of spacetime, most reasonably in the order $(t,x,y,z)$. Then the $t$-spinor is acted on by gravity but not the weak force, so is
a `right-handed' spinor in usual parlance. The $x$-, $y$- and $z$-spinors are then all `left-handed' spinors, acted on by both gravity and the weak force.
The allocation to space directions $x$, $y$ and $z$ is in some sense arbitrary, but should preferably be chosen with the ambient gravitational field in mind. 

The non-commuting of $SU(2)_G$ and $SU(2)_W$ implies that the gravitational eigenstates of particles are not compatible with the weak eigenstates.
Gravitational eigenstates are usually called `mass' eigenstates, and weak eigenstates are `flavour' eigenstates. In the case of neutrinos, the incompatibility
of `mass' and `flavour' eigenstates is well-known, and is the standard explanation for neutrino oscillations. Thus this model includes neutrino oscillations as
standard (unlike the Standard Model, which had to add them in later).

In the case of quarks, it is known that there is an incompatibility between the weak and strong eigenstates. It makes sense to regard the strong eigenstates as
mass eigenstates, so that we again see a discrepancy between mass eigenstates and flavour eigenstates. In the case of electrons, no such discrepancy
exists in the Standard Model, but it does in the proposed model. In this case, we must interpret the weak eigenstates as `inertial mass' eigenstates, since the
masses are measured by electromagnetic means, and interpret the gravitational eigenstates as `(active) gravitational mass' eigenstates.
Thus the model predicts that inertial mass and active gravitational mass are not the same concept. 

In other words, it predicts failure of the
(strong form of the) Weak Equivalence Principle. Note however that this does \emph{not} imply the inequivalence of inertial mass and
\emph{passive} gravitational mass. This weaker form of the Weak Equivalence Principle is called `universality of freefall', and is
strongly supported by experiment, but does not
contradict the proposed model. Unfortunately it is not always easy to tell in the literature which form of the Weak Equivalence Principle is being used in
any particular context. The reason for this appears to be an unquestioning acceptance of Newton's Third Law, which asserts the equivalence of active and passive 
gravitational mass for an \emph{instantaneous} gravitational action at a distance. However, if it is assumed, as is reasonable, that gravity travels at a finite speed,
most likely the same as the speed of light, then Newton's Third Law does not hold. 

Both $SU(2)_W$ and $SU(2)_G$ can be given the discrete structure of the binary tetrahedral group $Q_8\rtimes Z_3$, by taking generators
\begin{align}
\begin{pmatrix}1&0&0&0\cr 0&i&0&0\cr 0&0&j&0\cr 0&0&0&k\end{pmatrix}, \quad
\begin{pmatrix}1&0&0&0\cr 0&0&1&0\cr 0&0&0&1\cr 0&1&0&0\end{pmatrix}
\end{align}
for $SU(2)_W$, and
\begin{align}
\begin{pmatrix}i&0&0&0\cr 0&i&0&0\cr 0&0&i&0\cr 0&0&0&i\end{pmatrix},\quad
\begin{pmatrix}v&0&0&0\cr 0&0&v&0\cr 0&0&0&v\cr 0&v&0&0\end{pmatrix}
\end{align}
for $SU(2)_G$. In the Standard Model, the quaternions are reduced to complex numbers using the pure imaginary $t=v-w=i+j+k$.
In the case of $SU(2)_W$, this restriction has the effect of restricting to what is usually called the `third component' of weak isospin,
and ignoring the other two components. In the case of $SU(2)_G$, it has the effect of reducing to a scalar mass term, as defined by the
Dirac equation with a pure imaginary mass, again in $i+j+k$. In the notation of \cite{finite}, this corresponds to restricting from a
quaternionic mass defined by four particles (electron, proton, neutron and muon) to a mass defined by the neutron alone.

In other words, General Relativity is the effective theory that is obtained from this model in the continuous limit, under the assumption that all
matter consists of neutrons. The usual assumptions that the (active gravitational) mass of the proton is effectively equal to the neutron mass,
and that the electron mass is effectively zero, are only approximately true,
and may be invalid on scales much larger than the Solar System.
Moreover, the model implies that the muon mass cannot be ignored, since mass is a $4$-dimensional rather than $3$-dimensional concept.

\subsection{Experimental tests}
The model predicts, at a combinatorial level, a mixing between gravity and the other forces. In order to test this experimentally, however, it is first necessary to make
quantitative predictions. In a completely uniform gravitational field, the model reduces to the Standard Model, so the issue is to predict quantitatively
`anomalies' that should arise in a non-uniform gravitational field. Such anomalies can, of course, be expressed as
general relativistic corrections to the Standard Model, so the important question is, how large such corrections could be, and could they be detected experimentally?

Similarly, if all matter consists of neutrons, then the gravitational gauge group $SL(4,\HH)$
reduces to $SL(4,\RR)$, and the gravitational model reduces to General Relativity. But $SL(4,\RR)$ commutes with a copy of $\HH$ that allocates matter to four
essentially different types of particles. In \cite{finite,remarks} it is suggested that in addition to the `ordinary' particles (neutron, proton and electron) it is necessary to
add the muon to make up the four particles, and that a `mixing' between the electron and the muon interferes with the normal process of measuring mass.
The reason for this is the experimentally confirmed oscillation between electron neutrinos and muon neutrinos, which makes it impossible to tell, theoretically, whether
a given neutrino arose from an electron interaction or a muon interaction. 

First let us consider the effects of non-uniform gravitational fields on particle physics experiments. For experiments conducted in 
horizontal apparatus such as the Large Hadron Collider (LHC) we should not expect to see any measurable dependence on the strength of gravity,
but we might expect to see a dependence on the direction of gravity. From one side of the LHC to the other, this direction changes by nearly $4'$, but individual
experiments measure the effects of collisions in a relatively small part of the apparatus. 
On this scale,
the largest gravitational effects we could hope for are proportional to the sine of the angle by which the gravitational field changes direction. This is
roughly 1.6ppm for every 10 metres.

Anomalies of this order of magnitude have been reported in muon $g-2$ experiments \cite{muontheory,muonHVP,Fodor}, in which measurements are taken
over a distance of about $14$ metres. It is therefore possible to test whether a larger experiment would give a larger anomaly, as predicted by the 
gravitational model, or not. Another anomaly of the same order of magnitude is the CP-violation of neutral kaon decay \cite{CP}, measured over a
distance of approximately $17$ metres. In this case, a larger
experiment cannot distinguish between the gravitational model and more conventional models,
but an experiment oriented vertically rather than horizontally could do so.

Neutrino oscillation experiments are normally much larger, long-baseline, experiments, and the direction of the gravitational field changes
by a large angle, which may or may not be known, depending on the experiment. Testing for a dependence of experimental results on this angle
is something that could be done. Such tests could provide evidence for or against the hypothesis that neutrino oscillations are effected 
(or at least affected) by quantum gravity.

In the other direction, we are looking for effects of particle physics on (active) gravitational mass. Thus we are looking for systematic variations that
depend on the proportion of mass that is attributed to electrons. For this purpose we need to know how large a discrepancy between gravitational and
inertial mass of an electron could plausibly be, in experiments that we can do on Earth. 
This issue is discussed in detail in \cite{universal}, and appropriate experimental tests proposed. In particular, it is suggested that direct
measurements of the Newtonian gravitational constant $G$
are now accurate enough that slight differences between different materials, such as copper and gold, can be detected.
That is, if the inertial masses, and centres of mass, of a lump of copper and a lump of gold are measured accurately, and used to calibrate an 
experiment to
measure $G$, then a systematic difference between the two values, predicted by the model, should be large enough to be detected, 
somewhere in the range of 10 to 50 ppm,
depending on the precise assumptions made.

\subsection{Speculation}
The hope is, therefore, that by not breaking the symmetry of the quaternions, we can also avoid
the necessity to choose a definition of inertial frame. Or, to put it another way, the finite group implies that there is a `mixing' between the
generation symmetry and the choice of inertial frame. This type of mixing between particle physics and gravity was first considered by
Einstein \cite{Einstein1919} in 1919, but has not found favour in mainstream physics, in which the separation between gravity and
particle physics has been rigorously maintained. 
Nevertheless, the conventional choice of the laboratory frame of reference as `inertial' is fraught
with problems, particularly in the case of `long-baseline' 
experiments spanning two different laboratories,
with very different definitions of inertial frame. 

The model proposed here implies that such a change in
inertial frame should be reflected in a change in `generation' of fermions. Of course, there is no
experimental evidence that the generation of an electron is affected in this way, but there is some
experimental evidence that the generation (flavour) of a neutrino is \cite{SNO}. 
Indeed, it is shown in \cite{finite} that by inserting the masses of the three generations of electron
into the representations of the $Z_3$ quotient of our finite group, the experimental value of the
mixing angle between electron neutrino and muon neutrino can be obtained. In other words, the model implies that the
electron and neutrino generations are coupled together in a more complicated way than in the Standard Model.
 
As regards quarks, the model we are considering here allows us to put the masses of three generations of electron
into a quaternion representation of $Q_8\rtimes Z_3$, together with a proton or neutron, from which
the correct experimental value for the mixing angle between second and third generation quarks can be calculated \cite{finite}.
Again, this suggests a more complicated relationship between lepton and quark generations than is implied by the Standard Model.

There are in total a further $15$ parameters in the coordinates of the $4$-dimensional quaternionic representation on the spinors.
If we take at face value the proposed use of this representation for a theory of gravity, then we should expect these parameters to be
masses of elementary particles. 
Finally, note that by restricting the $4$-dimensional representation to \emph{subgroups} $Z_3$, $Q_8\rtimes Z_3$ and $G_{27}\rtimes Q_8\rtimes Z_3$
we obtain the first three representations again. It is therefore entirely possible that the $9$ mixing angles of the Standard Model are redundant
parameters, and can be derived from the $15$ masses. This is known to be the case for the electro-weak mixing angle, which is derived from the
mass ratio of the $W$ and $Z$ bosons. In \cite{finite} 
and Section~\ref{moreangles} we have discussed
six further possibilities.

\section{Conclusion}
\subsection{Summary}
The suggestion that physics might be fundamentally discrete, and that the apparently continuous behaviour we often observe is an illusion, goes back a long way.
In the context of quantum mechanics, this suggestion was made by Einstein \cite{Einstein1935} as early as 1935, but no discrete theory was found to rival the
continuous models based on the Schr\"odinger equation and then the Dirac equation. The Standard Model of Particle Physics, largely completed by the mid 1970s,
is entirely based on continuous groups, so that the emergence of discrete experimental observations is still not explained in a totally convincing way \cite{Bell,Kochen,Rae}.

Theories beyond the Standard Model are almost without exception continuous theories. 
A few lone voices, such as 't Hooft \cite{thooft,thooftbook}, still speak out in favour of
Einstein's dream of a discrete foundation, even if this has to be at the Planck scale.
But a convincing discrete model is still lacking, and few people are seriously looking for one. 
The ones I have put forward in \cite{finite,octahedral,gl23model,icosa,tetrions,ZYM}, using either the binary tetrahedral group,
the binary icosahedral group or one of the two isoclinic variants of the binary octahedral group, all suffer
from drawbacks of various kinds, most notably the lack of any convincing analogue of the colour $SU(3)$,
which is an essential part of the Standard Model.
Nevertheless, it is a mathematical fact that finite groups are \emph{much} more complicated than Lie groups,
and quite small finite groups can model properties that cannot be modelled with Lie groups at all.

In this paper, therefore, I have tried to remedy this defect by explicitly proposing a finite analogue or precursor of $SU(3)$.
At a combinatorial (qualitative) level it has many desirable properties that mirror known experimental properties of elementary
particles, in particular the chirality of the weak interaction, the mixing of electromagnetism and the weak and strong nuclear forces,
and neutrino oscillations. 
Quantitative analysis of masses, mixing angles and probability amplitudes has only been touched upon,
and is in general 
 beyond the scope of this paper, as it requires detailed investigations
of coordinates for the representations. 

\subsection{Further work}
The proposed unification of Gell-Mann and Pauli matrices using a 3-dimensional complex reflection group really
does little more than provide a discrete foundation for the Standard Model of Particle Physics as it currently stands.
Nevertheless, it permits new methods of calculation, using finite geometries, which have been shown to be effective in calculating six or 
(probably) seven of the
mixing angles from fundamental mass parameters. 
A priority for further work is therefore to try to find similar calculations for the remaining two mixing angles.
It is possible that this cannot be done with the Pauli and Gell-Mann matrices alone, but
requires the Dirac matrices as well.

Such unification with the Dirac matrices 
cannot be done in the Standard Model, since it is impossible
to arrange for the Lorentz group to commute with the whole of the gauge group in this case.
Nevertheless, the non-commutation is very subtle, and only really affects properties of neutrinos and gluons, as well
as the Higgs field, so that any differences from a commuting theory are 
hard to detect.
What this means is that a unification with Dirac matrices provides a possible mechanism for particle physics to influence the shape of spacetime,
and hence couple to gravity. In particular, the proposed unification using a quaternionic reflection group provides a potential foundation
for a quantum theory of gravity. 
Investigating this possibility further is a priority for future work. 

The appropriate compact gauge group here is $Sp(4)$,
of dimension $36$, but it is also possible to gauge the theory with the noncompact group $SL(4,\HH)$ of dimension $63$. In this case,
there is a subgroup $SL(4,\RR)$ that might be used for general covariance, in order to reduce any such putative quantum theory of gravity
to General Relativity in the appropriate limit. In the context of the Coleman--Mandula Theorem, therefore, the gauge group
for quantum gravity reduces to $SU(2)$, that is the centralizer of $SL(4,\RR)$ in $SL(4,\HH)$, as has also been proposed in other models, such as  \cite{WoitSO24}.

Since $SL(4,\RR)$ contains no quaternions, and does not contain this gauge group $SU(2)$,
it cannot model the three generations of electrons.
The model therefore suggests that it is the (quantum and/or classical)
gravitational properties of
electrons \cite{universal} that need to be investigated in case of any apparent deviations \cite{Rubin,Chae} of gravity from General Relativity.

If a quantum theory of gravity along these lines turns out to be viable, then it might lead to methods for calculating the $15$ mass parameters
of the Standard Model from something more fundamental. It is not at all clear what could be `more fundamental', or how it could be used,
so that this is at present a highly speculative idea. Nevertheless, it is clear that it will be necessary to examine closely the theoretical and practical definitions
of inertial frames of reference in different contexts, and to calculate the effects of transformations between them. 
Some speculations in this direction are presented in \cite{icosa,universal}, in which a preliminary 
attempt is made to relate the laboratory `inertial' frame to the Solar System `inertial' frame.

Since $SL(4,\HH)$ is a real form of $SU(8)$, of Lie type $A_7$, and since $A_7$ is contained in both $E_7$ and $D_8$, 
there are at least two ways to embed this gauge group in (complex) $E_8$. It may even be possible to embed it in a real form of $E_8$.
In the $E_7$ case, there are no spinor representations of $SL(4,\HH)$ in $E_8$, which seems to rule this one out \cite{QE7,DMWE7}.
In the $D_8$ case, the spinors split into $4+28$ quaternionic dimensions, which may be viable as a splitting into neutrinos versus the charged
fundamental fermions. Even so, the model would 
still be quite unlike the various $E_8$ models \cite{octions,Chester,Lisi,Lisi2} that have been proposed over a number of years.

\end{document}